\documentclass[11pt]{amsart}

\usepackage{amsmath,amsthm,amssymb,mathrsfs,amsfonts,verbatim,enumitem,color,leftidx}
\usepackage{mathabx}
\usepackage{etoolbox} 

\usepackage{bbm}
\usepackage[all,tips]{xy}
\usepackage{graphicx,ifpdf}
\usepackage{stmaryrd}
\usepackage[a4paper,left=20mm,right=20mm,top=15mm,bottom=25mm]{geometry}
\ifpdf
\fi
\usepackage[colorlinks]{hyperref}
\hypersetup{
linkcolor=blue,          
citecolor=green,        
}

\usepackage{tikz}

\usepackage{marginnote}
\usepackage{color}

\newtheorem{thm}{Theorem}
\newtheorem{lem}{Lemma}
\newtheorem{coro}{Corollary}
\newtheorem{prop}{Proposition}

\newtheorem{rmk}{Remark}

\numberwithin{equation}{section}

\theoremstyle{remark}
\numberwithin{equation}{section}

\definecolor{esperance}{rgb}{0.0,0.5,0.0}

\newcommand{\onto}{\xymatrix{\ar@{>>}[r]&}}

\usepackage{amssymb}
\usepackage{amsmath,amssymb,graphicx}
\usepackage[utf8]{inputenc}
\usepackage{subcaption}

\usepackage{pgfplots}
\pgfplotsset{compat=1.14}
\pgfplotsset{my style/.append style={axis x line=middle, axis y line=
middle, xlabel={$\frac{1}{p}$}, ylabel={$\alpha$}, axis equal }}
\usepackage[utf8]{inputenc}

\title{Newton polyhedrons and $L^p$ Sobolev estimations}
\author{Kiseok, Yeon }

\begin{document}
\maketitle
\begin{abstract}
    The aim of this study is to provide a perspective to help understand the singular average operator over polynomial hypersurfaces. In particular, this perspective will provide brevity and the possibility of generalizing previous results dealing with the fundamental problem of determining the precise $L^p$ regularity enhancement for the average operators. In previous studies dealing with polynomials, the Newton polyhedron of a polynomial has been utilized to observe dominant monomials. In this study, we go further by discussing the involvements of other monomials in detail, by introducing several geometric values on the Newton polyhedron.
\end{abstract}
\section{Introduction}
 This study considers the fundamental problem of determining the precise $L^p$ regularity enhancement under convolution, using a singular measure over certain polynomial hypersurfaces in $R^3$.
 In [3], polynomials of the form $t_1^2+t_1^at_2^b+ t_2^M$ ($a,b,M$ positive even numbers) are considered, and the optimal $L^p \rightarrow L^p_{\alpha}$ ranges of their operators are represented as ($\frac{1}{p},\alpha(p))$ graphs according to the arrangement of $(a,b)$, $(2,0)$,  and $(0,M)$. In this study, our aim is to investigate how several geometric values on the Newton polyhedron explain this arrangement, and to obtain optimal $L^p\rightarrow L^p_{\alpha}$ ranges for more general integral operators. However, unlike in [3] our focus is not on the endline estimates. 
 
 Let us define an operator $\mathcal{A}$ by
$$ \mathcal{A}(f)(x_1,x_2,x_3)=\int f(x_1-t_1,x_2-t_2,x_3-P(t_1,t_1))\psi(t_1,t_2)dt_1dt_2,$$

 where $\psi(t_1,t_2)$ is a $C^{\infty}$ function with sufficiently small support near the origin, and  
$P(t_1,t_2)=\displaystyle\sum_{(m,n)}a_{mn}t_1^mt_2^n$ ($m,n \neq1 $) satisfies certain conditions, which will be introduced in (2.1) and (2.2).
\section{Statement of main result}

Let $\mathcal{P}(\mathbb{R}^2)$ be the family of all polynomials on $\mathbb{R}^2$. For $P\in\mathcal{P}(\mathbb{R}^2)$, the Newton polyhedron $N(P)$ is the smallest convex set (convex hull) containing the union of the first quadrants translated by the exponents $\Lambda(P)$ of the all monomials in $P$. For  $P(t_1,t_2)=\displaystyle\displaystyle\sum_{(m,n)\in \Lambda(P)}a_{mn}t_1^mt_2^n$ with $\Lambda(P)=\{(m,n): a_{mn} \neq 0\}$, we write 

$\quad$
$$N(P)=CH( \bigcup_{(m,n)\in {\Lambda(P)}} (m,n)+{(X,Y): X\geq 0, Y\geq 0}).$$
For $\vec{v}$ in $\mathbb{R}^2_+$, we define $N_{\vec{v}}(P)$ using $\mathbb{V}^{\bot}=\{\vec{w}\in \mathbb{R}^2 :<\vec{v},\vec{w}>=0\}$:

$$N_{\vec{v}}(P)=CH(\bigcup_{(m,n)\in \Lambda(P)}(m,n)+ \mathbb{V}^{\bot}).$$

The Newton polygon $N(P)$ is an unbounded convex polygon, with finite numbers of edges and vertices. By $\mathcal{E}(P)$ and $\mathcal{V}(P)$, we denote the families of edges and vertices of $N(P)$, respectively, and we denote the boundary of $\mathcal{N}(P)$ by $F(P)$. For each edge $E\in \mathcal{E}(P)$, we let $P_E$ be the polynomial whose exponents of nonzero monomials of $P$ are restricted to be on $E$:
$$P_E(t_1,t_2)= \displaystyle\displaystyle\sum_{(m,n)\in E\cap \Lambda(P)}a_{mn}t_1^mt_2^n.$$

For the reason that our focus is on figuring out the relationship between $L^p$  regularity enhancement of the operator $\mathcal{A}$ and $\Lambda(P)$, we impose the following two conditions related to coefficients $a_{mn}$'s. In fact, two conditions are given on $P_E(t_1,t_2)$'s. To resolve this problem without any conditions, we will have to study $L^p$ regularity enhancement of the operator $\mathcal{A}$ corresponds to $P_E(t_1,t_2)$.

$\ $

$\bullet$ Let $E_{t_1},E_{t_2}$ be the edges containing only monomials whose exponents are of the form $(m,0), (0,n)$, respectively. By  $\mathcal{E}^{\circ}(P)$, we denote the set of edges in $\mathcal{E}(P)$ excluding $E_{t_1}, E_{t_2}.$

\ 
$\ $


For all $E\in \mathcal{E}^{\circ}(P), E_{t_1}, E_{t_2}$, 
\begin{equation}
| det(\nabla^2 P_E(t_1,t_2))|>0,|\partial_{t_1}^2P_{E_{t_1}}|,|\partial_{t_2}^2P_{E_{t_2}}|>0\ \textrm{for}\ |(t_1,t_2)|\neq0\
\end{equation}

$\bullet$ Let $P^1,P^2$ be polynomials after removing monomials whose exponents lie on $E_{t_1},E_{t_2}$ from $P(t_1,t_2)$, respectively. 
$$\Lambda(P^1)=\{(m,n)\in\Lambda(P)|n>0\}$$
$$\Lambda(P^2)=\{(m,n)\in\Lambda(P)|m>0\}$$

 
 For all $ E_1\in \mathcal{E}(P^1)\setminus\mathcal{E}(P)$ and $
E_2\in\mathcal{E}(P^2)\setminus\mathcal{E}(P) $,
\begin{equation}
|\partial_{t_2}^2P_{E_1}(t_1,t_2)|>0, |\partial_{t_1}^2P_{E_2}(t_1,t_2)|>0\  \textrm{for } |(t_1,t_2)|\neq0 \
 \end{equation}

$\bullet$ Examples
 
  $(1)\ \textrm{Let}\ P(t_1,t_2)= t_1^2+t_1^at_2^b+t_2^M $ ($a,b,M$ positive even numbers) satisfies conditions 1 and 2. Indeed, as $P_E=t_1^2+t_2^M$, $P_{E_{t_1}}=t_1^2$, $P_{E_{t_2}}=t_2^M$, and $P^1=t_1^at_2^b+t_2^N$, we obtain $| det(\nabla^2 P_E(t_1,t_2))|=|2M(M-1)t_2^{M-2}|>0,|\partial_{t_1}^2P_{E_{t_1}}|=2>0,|\partial_{t_2}^2P_{E_{t_2}}|=|M(M-1)t_2^{M-2}|>0 \ \textrm{for}\ |(t_1,t_2)|\neq0$. Furthermore, as $P_{E_1}=t_1^at_2^b+t_2^M$ and $P_{E_2}= t_1^at_2^b+t_1^2$, we obtain $|\partial_{t_2}^2P_{E_1}(t_1,t_2)|=|b(b-1)t_1^at_2^{b-2}+M(M-1)t_2^{M-2}|>0, |\partial_{t_1}^2P_{E_2}(t_1,t_2)|=|a(a-1)t_1^{a-2}t_2^b+2|>0\  \textrm{for } |(t_1,t_2)|\neq0$.

$(2)$ $\ \textrm{Let}\ P(t_1,t_2)=t_1^{m_0}+\displaystyle\sum_{(m,n)}a_{mn}t_1^mt_2^n$ ($m_0,m,n$ positive even numbers, $a_{mn}>0$). Assume that $m_0$ is the smallest value among the exponents of $t_1$. Then, $P(t_1,t_2)$ satisfies conditions 1 and 2. Indeed, as all exponents are positive even integers, the second derivative of each of monomial in $P(t_1,t_2)$ is positive, as in example (1), for $|(t_1,t_2)|\neq 0$. Furthermore, as $\mathcal{E}^{\circ}(P)$ is the empty set, $P(t_1,t_2)$ satisfies condition 1.

$\ $

Viewing a polynomial in light of the Newton polyhedron, we can sort the dominant monomials by size according to the sections of the domain over which our integrals are defined. In this study, we show that not only are dominant monomials involved in $L^p$ regularity, but second dominant monomials also affect the $L^p$ regularity of our operator. Furthermore, the Newton polyhedrons of these second dominant monomials play a significant role in sharply determining the $L^p\rightarrow L^p_{\alpha}$ ranges. To demonstrate this involvement, we define several geometric values on the Newton polyhedron.

$\ $

$\mathbf{Definition\ 1.}$ Let $\vec{v}$ be a vector in $\mathbb{R}_+^2$. The boundary of the polygon $\mathcal{N}_{\vec{v}}(P)$ and $\mathcal{N}(P)$ intersect with the diagonal line $\{(X,Y): X=Y\}$ at the unique points $(\delta_{\vec{v}}(P),\delta_{\vec{v}}(P))$ and $(\delta(P),\delta(P))$, respectively. The values $\delta_{\vec{v}}(P)$ and $\delta(P)$ are defined by
$$\delta_{\vec{v}}(P)= \textrm{min}\{t: (t,t)\in N_{\vec{v}}(P)\},$$
$$\delta(P)= \textrm{min}\{t: (t,t)\in N(P)\}.$$

We denote $\delta(P)$, called the Newton distance for $P(t_1,t_2)$, by $\delta$, to simplify the notation.

$\mathbf{Definition\ 2.}$ Define $m_s$ and $n_s$ by
$$m_s= \displaystyle\min_{(m,0)\in \Lambda(P)}m,$$ $$n_s=\displaystyle\min_{(0,n)\in \Lambda(P)}n.$$ If $\Lambda(P)$ contains no elements of the form of $(m,0)$ (respectively $(0,n)$), then we take $m_s=\infty$ (respectively $n_s=\infty$). 

$\mathbf{Definition\ 3.}$ Let $M$ be $\displaystyle\min_{(m,n)\in\Lambda(P^1)}m$. By $(\delta^{(m,n)},\delta^{(m,n)})$, denote the intersection point between $y=x$ and the ray $\overrightarrow{(m_s,0),(m,n)}$. Define $(M_s,N_s)$ by the point closest to $(m_s,0)$ among those points satisfying  
$$\delta^{(M_s,N_s)}\leq \delta^{(m,n)}\ \textrm{for all}\ (m,n)\in  [F(P^1)\cap\Lambda(P^1)]\cup(M,\infty).$$
(See Figure 1 for examples.)
\begin{figure*}
  \includegraphics[width=0.49\textwidth]{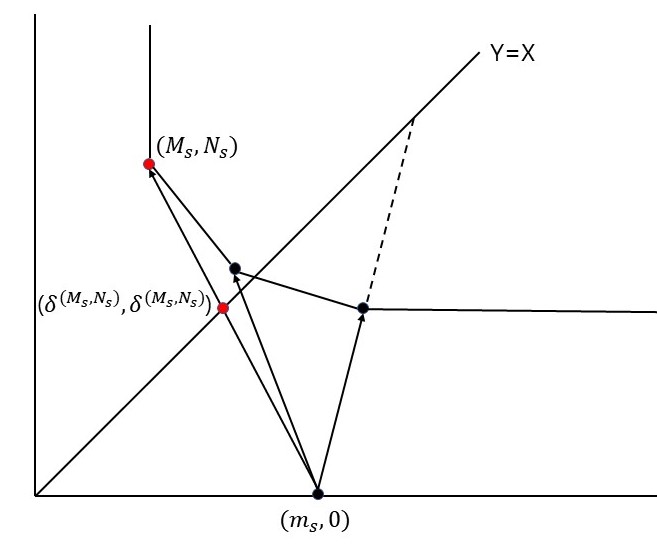}
  \label{1}
    \includegraphics[width=0.49\textwidth]{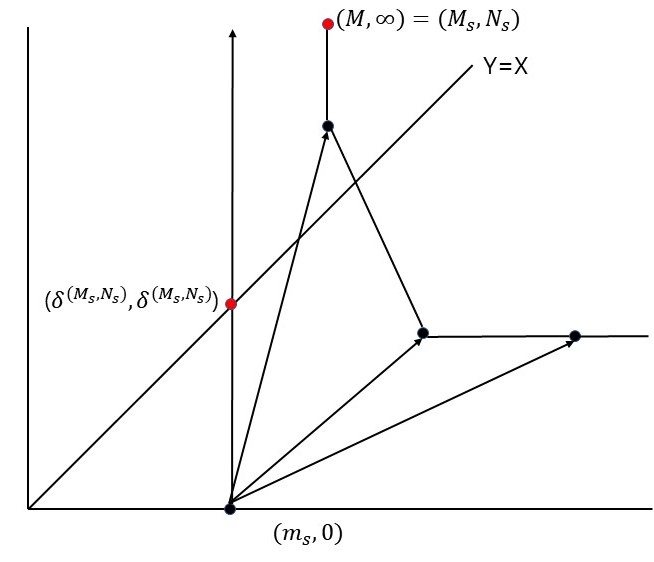}
    \label{9}
\caption{Examples for $\delta^{(M_s,N_s)}$}  
\end{figure*}
$\ $

$\mathbf{Definition\ 4.}$ For $(M_s,N_s)$ and $m_s$, define a set $\mathcal{R}$ by
$$\mathcal{R}=\{(m,n)\in [F(P^1)\cap\Lambda(P^1)]\cup(M,\infty) |  2\delta^{(m,n)}\leq n\leq N_s\}.$$

(See Figure 2 for examples.)

$\ $

In the case that $\mathcal{R}\neq\phi$, we set $(M_s,N_s)=(M_1,N_1)$, and enumerate the vertices in $\mathcal{V}(P^1)$ as 
$(M_2,N_2),\cdots (M_{n+1},N_{n+1})$, so that
$(M_k,N_k)\in \mathcal{R}$ for $k=2,..,n$, 
$ (M_{n+1},N_{n+1})\not\in\mathcal{R}$, and $N_1\geq N_2\geq\cdots\geq N_{n+1}$.
By $\vec{v}_i$, we denote a normal vector of 
$ \overline{(M_i,N_i)(M_{i+1},N_{i+1})}$ in $\mathbb{Z}_+^2$.

\begin{figure*}
  \includegraphics[width=0.8\textwidth]{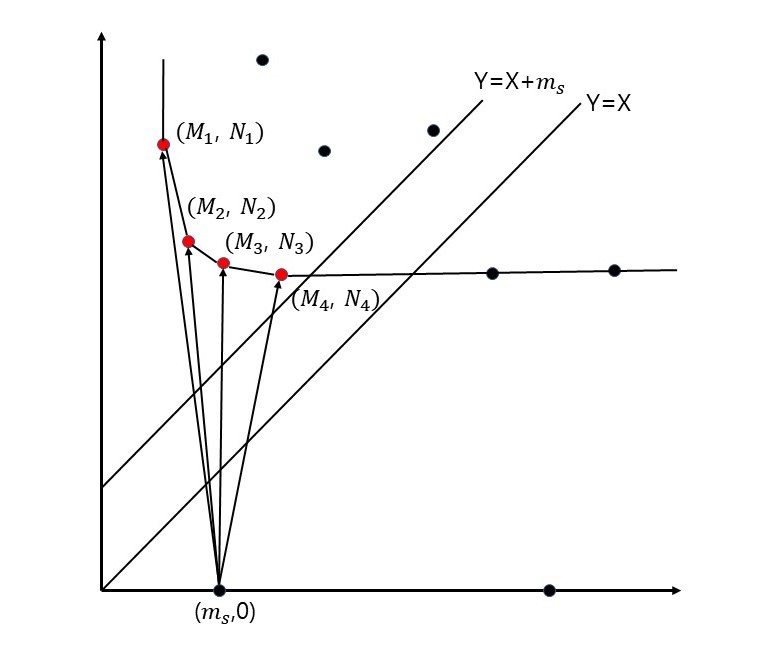}
\caption{$\mathcal{R}$} 
\label{2}
\end{figure*}
\begin{rmk}

(1) We may assume that our polynomial satisfies 
$$m_s = \min(m_s,n_s),$$ as otherwise we may exchange $t_1$ and $t_2$.

$\ $

(2) We can observe that $\delta^{(M_1,N_1)}<\delta^{(M_2,N_2)}<\cdots<\delta^{(M_n,N_n)}$.


\end{rmk}
$$\ $$
We classify all possible arrangements of exponents of the polynomial into two situations by considering $\mathcal{R}$. Specifically, $\mathcal{R}=\phi$ and $\mathcal{R}\neq\phi$.


\begin{thm}

 Assume that a polynomial $P(t_1, t_2)$ satisfies (2.1) and (2.2). Then, the operator $\mathcal{A}$ maps $L^p$ to $L^p_{\alpha}\  $for $\alpha<\alpha(p)$ under the following conditions.

$\ $

$\ $

$\bullet$ In the case that $ \mathcal{R}=\phi$:

$\ $

$\ $

$\alpha(p):=\left\{\begin{array}{l} \frac{1}{\delta}\quad\qquad\textrm{if}\  \frac{1}{2\delta}<\frac{1}{p}<\frac{1}{2}\\ 
\frac{2}{p} \quad\qquad \textrm{if}\ \frac{1}{p} \leq \frac{1}{2\delta}\end{array}\right.$

$$\ $$
$\ $

$\bullet$ In the case that $\mathcal{R}\neq \phi$:

$\ $

$\alpha(p):=\left\{\begin{array}{l} \frac{1}{\delta}\quad\qquad\qquad\textrm{if}\  \min(\frac{1}{\delta^{(M_s,N_s)}}-\frac{1}{N_s},\frac{1}{2})\leq\frac{1}{p}\leq\frac{1}{2}\\ l^{\vec{v_1}}(\frac{1}{p})\qquad\quad\textrm{if}\ \frac{1}{\delta^{(M_2,N_2)}}-\frac{1}{N_2}\leq\frac{1}{p}\leq\min(\frac{1}{\delta^{(M_s,N_s)}}-\frac{1}{N_s},\frac{1}{2})\\
\cdot\qquad\qquad\qquad\qquad\qquad\cdot\\
\cdot\qquad\qquad\qquad\qquad\qquad\cdot\\
l^{\vec{v}_{n-1}}(\frac{1}{p})\  \qquad\textrm{if}\ \frac{1}{\delta^{(M_n,N_n)}}-\frac{1}{N_n}\leq\frac{1}{p}\leq \frac{1}{\delta^{(M_{n-1},N_{n-1})}}-\frac{1}{N_{n-1}}\\
l^{\vec{v}_n}(\frac{1}{p})\ \qquad\quad\textrm{if} \frac{1}{\delta_{\vec{v}_n}(t_1^{m_s})+\delta_{\vec{v}_n}(t_1^{M_n}t_2^{N_n})}\leq\frac{1}{p}\leq \frac{1}{\delta^{(M_n,N_n)}}-\frac{1}{N_n}\\
\frac{2}{p}\  \quad\qquad \qquad\textrm{if}\ \frac{1}{p} \leq \frac{1}{\delta_{\vec{v}_n}(t_1^{m_s})+\delta_{\vec{v}_n}(t_1^{M_n}t_2^{N_n})}\end{array}\right.$

where $l^{\vec{v_k}}(\frac{1}{p})=(1-\frac{\delta_{\vec{v}_k}(t_1^{m_s})}{\delta_{\vec{v}_k}(t_1^{M_k}t_2^{N_k})})\frac{1}{p}+ \frac{1}{\delta_{\vec{v}_k}(t_1^{M_k}t_2^{N_k})}\ (k=1,2,\cdot\cdot\cdot,n)$.
\end{thm}

$$\ $$

\begin{rmk}
(1) The indicated ranges of the parameters $p$ and $\alpha$ in Theorem 1.1 cannot be improved, in the sense of the unboundness of the $L^p \rightarrow L^q$ estimates under affine transformations between the $L^p\rightarrow L^p_{\alpha}$ and $L^p \rightarrow L^q$ bounds.

(2) We cannot determine whether $\mathcal{A}$
satisfies $L^p\rightarrow L^p_{\alpha(p)}$ at some of the endlines.

(3) It holds that $l^{\vec{v_k}}(\frac{1}{\delta^{(M_k,N_k)}}-\frac{1}{N_k})=\frac{1}{\delta^{(M_k,N_k)}}$ for k=1,2,$\cdot\cdot\cdot$,n, which will be proved in Lemma3.
\end{rmk}

$$\ $$

$\mathbf{Sketch\  of\ the\  proof\ of\ Theorem\ 1}$

$\ $

$\bullet$\ We set $\Lambda(P)= \{(m,0),(M,N)\}$ ($m \neq 1,M,N \geq 2$), so that the polynomial $P$ satisfies $(2.1)$ and $(2.2).$

For two fixed vectors $\vec{v},\vec{w}(\in \mathbb{Z}_+^2)$, consider the operators $$\mathcal{A}_I(f)(x)=\int f(x_1-t_1,x_2-t_2,x_3-P(t_1,t_2))\displaystyle\sum_{(j_1,j_2)\in I}\eta(\frac{t_1}{2^{-i_1}},\frac{t_2}{2^{-i_2}})dt$$
$$\mathcal{A}_J(f)(x)=\int f(x_1-t_1,x_2-t_2,x_3-P(t_1,t_2))\displaystyle\sum_{(j_1,j_2)\in J}\eta(\frac{t_1}{2^{-j_1}},\frac{t_2}{2^{-j_2}})dt,$$ where
$I=\{\vec{v}i|i\in\mathbb{Z}_+\}$, $J=\{\vec{v}i+\vec{w}l|i,l\in\mathbb{Z}_+\}$, and $\textrm{supp}(\eta(t_1,t_2))$ is contained in $[\frac{1}{2},1]\times [\frac{1}{2},1]$.

$\ $

$\ $

$\bullet$ Let $(\vec{n}_k)_{k}$ be the collection of the normal vectors of $\mathcal{E}(P)$ and $\mathcal{E}(P^1)$ in $\mathbb{Z}_+^2$. We arrange $\vec{n_k}(k\geq1)$ in ascending order by their slopes. Using these vectors $\vec{n}_k$, we decompose $\mathbb{Z}_+^2$ as $$\mathbb{Z}^2_+=(\displaystyle\bigcup_{k}\{\vec{n}_ki|i\in \mathbb{Z}_+\})\cup(\displaystyle\bigcup_{k}\displaystyle\bigcup_{\vec{s}\in \mathcal{S}}\{\vec{s}+\vec{n}_ki+\vec{n}_{k+1}l| i,l\in \mathbb{Z}_+\}),$$ 
 where $\vec{s}$ is a vector in $\mathbb{Z}^2$. Furthermore, $|\mathcal{S}|$ is finite and depends on $\vec{v}$ and $\vec{w}$. Then, we obtain the $L^p\rightarrow L^p_{\alpha}$ ranges for our operator $\mathcal{A}$ by taking the intersection of the $p$ and $\alpha$ ranges for $\mathcal{A}_I$ and $\mathcal{A}_J$ for the decomposed index sets above.

$\ $

$\bullet$ Obtain necessary conditions for the $L^p\rightarrow L^p_{\alpha}$ ranges through the unboundness of the $L^p\rightarrow L^q$ estimates.

$$\ $$

$\mathbf{Notation}.$ We utilize $A \lesssim B$ if $A\leq cB$ for some constant $c$, and $A\sim B$ if $c_1B \leq A \leq c_2B$ for some constants $c_1, c_2>0$. We denote $(1,1)$ by $\vec{\mathbf{1}}$.

$\ $

In the following, we may assume that
$1\lesssim|\xi|,|\xi_1|+|\xi_2|\lesssim |\xi_3|$, and so $|\xi|\sim|\xi_3|$. This is because if $|\xi_1|+|\xi_2|>C|\xi_3|$ for some large constant $C$ or $|\xi|\lesssim 1$, then we have by integration by parts that

$|m(\xi)|\leq C_N(1+|\xi|)^{-N}$ for any $N$, where $\widehat{\mathcal{A}(f)}(\xi)=m(\xi)\hat{f}(\xi)$. 

Thus, we can set $$\|\mathcal{A}(f)\|_{L^p_{\alpha}}\sim \|[|\xi_3|^{\alpha}\widehat{\mathcal{A}(f)}]^{\vee}\|_{L^p}.$$ Let $\psi_0(t)$ be a $C^{\infty}_0(\mathbb{R})$ function such that $\psi_0(t)=1$ for $|t|\leq \frac{1}{2}$ and $\psi_0(t)=0$ for $|t|>1$. Define another function $\eta(t)$ by $\eta(t)=\psi_0(t)-\psi_0(2t).$ Then, 
$\displaystyle\sum_{j=1}^{\infty}\eta(2^jt)=1$, for all $t\leq\frac{1}{4}$.
$$\mathcal{A}(f)(x_1,x_2,x_3)= \int f(x_1-t_1, x_2-t_2, x_3- P(t_1,t_2))\psi(t_1,t_2)dt_1dt_2.$$
\begin{equation*}
 \begin{split}
\widehat{A(f)}(\xi)&= \hat{f}(\xi)\int e^{i(t_1\xi_1+t_2\xi_2+P(t_1,t_2)\xi_3)}\psi(t_1,t_2)dt \\
 &= \hat{f} (\xi)\displaystyle\displaystyle\sum\limits_{j_1,j_2=1}^{\infty}2^{-j_1}2^{-j_2}\int e^{i\xi_3\phi_J(t_1,t_2,\xi)}\eta(t_1,t_2)dt_1dt_2\\
 &= \hat{f}(\xi) \displaystyle\displaystyle\sum\limits_{j_1,j_2=1}^\infty 2^{-j_1}2^{-j_2} m_{j_1,j_2}(\xi_1,\xi_2,\xi_3),
 \end{split}
 \end{equation*}
 
where $\phi_J(t_1,t_2,\xi)=2^{-j_1}t_1\frac{\xi_1}{\xi_3}+2^{-j_2}t_2\frac{\xi_2}{\xi_3}+P(2^{-j_1}t_1,2^{-j_2}t_2)$ and $ \eta(t_1,t_2)=\psi(2^{-j_1}t_1, 2^{-j_2}t_2)\eta(t_1)\eta(t_2)$. We do not care about the index $j_1, j_2$ on $\eta(t_1,t_2)$ when we use the vander corput lemma. Indeed, it holds that $|\nabla \eta(t_1,t_2)|= |\nabla(\psi(2^{-j_1}t_1,2^{-j_2}t_2))\eta(t_1)\eta(t_2)+\psi(2^{-j_1}t_1,2^{-j_2}t_2)\nabla(\eta(t_1)\eta(t_2))|\sim 1$.

\section{Propositions and lemmas}

\begin{lem} 
For a fixed $\vec{n}=(n_1,n_2)$ in $\mathbb{Z}_+^2$, assume that $\delta_{\vec{n}}(t_1^{\tilde{m}}t_2^{\tilde{n}})<\delta_{\vec{n}}(t_1^mt_2^n)$ for all $(m,n)\in \Lambda(P)\backslash (\tilde{m},\tilde{n})$. Then, if $t_1\sim2^{-n_1i}$ and $t_2\sim2^{-n_2i}$ for a sufficiently large integer $i$, it holds that

$$\frac{1}{2}a_{\tilde{m}\tilde{n}}t_1^{\tilde{m}}t_2^{\tilde{n}}\leq\sum_{(m,n)\in\Lambda(P)}a_{mn}t_1^mt_2^n\leq 2a_{\tilde{m}\tilde{n}}t_1^{\tilde{m}}t_2^{\tilde{n}}.$$
\end{lem}
$\ $

$Proof$  By the definitions of $\delta_{\vec{n}}(t_1^{\tilde{m}}t_2^{\tilde{n}}))$ and $\delta_{\vec{n}}(t_1^mt_2^n)$, it follows that
$$\vec{n}\cdot((\tilde{m},\tilde{n})-\delta_{\vec{n}}(t_1^{\tilde{m}}t_2^{\tilde{n}}))\vec{\mathbf{1}})=0$$
$$\vec{n}\cdot((m,n)-\delta_{\vec{n}}(t_1^{m}t_2^n))\vec{\mathbf{1}})=0$$
Thus, $\vec{n}\cdot(m,n)-\vec{n}\cdot(\tilde{m},\tilde{n})=(\delta_{\vec{n}}(t_1^mt_2^n)-\delta_{\vec{n}}(t_1^{\tilde{m}}t_2^{\tilde{n}}))\vec{n}\cdot\vec{\mathbf{1}} >0$. Because there are finitely many $(m,n)\in \Lambda(P)\backslash(\tilde{m}, \tilde{n})$, there exists $\epsilon>0$ such that 

$$\epsilon =\min\{\delta_{\vec{n}}(t_1^mt_2^n)-\delta_{\vec{n}}(t_1^{\tilde{m}}t_2^{\tilde{n}})|(m,n)\in\Lambda(P)\backslash(\tilde{m}, \tilde{n})\},$$ which implies that there exists $\epsilon>0$ such that 
\begin{equation*}
\vec{n}\cdot (\tilde{m},\tilde{n}) +\epsilon <  \vec{n}\cdot(m,n) \ \textrm{for\ all\ }(m,n)\in\Lambda(P)\backslash(\tilde{m}, \tilde{n}).\qquad\qquad 
\end{equation*}
Thus, 
\begin{equation*}
\displaystyle\sum_{(m,n)\in\Lambda(P)}a_{mn}t_1^mt_2^n\sim\displaystyle\sum_{(m,n)\in\Lambda(P)}a_{mn}2^{-(n_1,n_2)\cdot(m,n)i}\sim a_{\tilde{m}\tilde{n}}2^{-(n_1,n_2)\cdot(\tilde{m},\tilde{n})i}
\sim a_{\tilde{m}\tilde{n}}t_1^{\tilde{m}}t_2^{\tilde{n}},
\end{equation*}
and the proof is complete.

$\ $

We call such a $t_1^{\tilde{m}}t_2^{\tilde{n}}$ the first dominant monomial for $\vec{n}$. We can then determine the second dominant monomial for ${\vec{n}}$ inductively.

$\ $
\begin{lem}
Let $P(t_1,t_2)=\displaystyle\sum a_{mn}t_1^mt_2^n$. Let $\vec{v}, \vec{w}$ be vectors in $\mathbb{Z}_+^2$. Define 
\begin{equation*}
    m_{i,l}(\xi)=\int e^{i\xi_3\phi_J(t_1,t_2,\xi)}\eta(t_1,t_2)dt_1dt_2,
\end{equation*}
where $ J=\{(j_1,j_2)=\vec{v}i+\vec{w}l|i,l\in \mathbb{Z}_+\}.$ and $supp(\eta(t_1,t_2))\subset([-1,-\frac{1}{2}]\cup[\frac{1}{2},1])\times([-1,-\frac{1}{2}]\cup[\frac{1}{2},1])$.
Then, 
$$\|(\sum_{i}|m_{i,l}^{\vee}\ast f_{i,l}|^2)^{\frac{1}{2}}\|_{p}\leq C\|(\sum_{i}|f_{i,l}|^2)^{\frac{1}{2}}\|_p,\ \textrm{and}\ C\ \textrm{does not depend on j}.$$
\end{lem}

$Proof$ It suffices to show that 
$$\|\displaystyle\sup_im_{i,l}^{\vee}\ast f\|_p\leq C|f\|_p\ \ \textrm{for}\ p>1\ \textrm{and C does not depend on}\ j.$$

Let $\tilde{P}(t_1,t_2)$ be a polynomial, which will be determined later.
Define $$\psi_{i,l}(\xi)=\int e^{i\xi_3\tilde{\phi}_J(t_1,t_2,\xi)}\eta(t_1,t_2)dt_1dt_2,$$ where $\tilde{\phi}_{J}(t_1,t_2,\xi)=2^{-j_1}t_1\frac{\xi_1}{\xi_3}+2^{-j_2}t_2\frac{\xi_2}{\xi_3}+\tilde{P}(2^{-j_1}t_1,2^{-j_2}t_2)$.

Let us define a Littlewood--Paley projection
\begin{equation*}
    \widehat{(P_{i,l,k}\ast f)}(\xi)=\chi(\frac{\xi_3}{2^{(m,n)\cdot(\vec{v}i+\vec{w}l)+k}})\hat{f}(\xi),
\end{equation*}
 where $\chi$ is supported on $[-1,-\frac{1}{2}]\cup[\frac{1}{2},1]$ and $\displaystyle\sum_{k}\chi^2(2^{-k}\xi_3)=1$. Then, 
\begin{multline*}
\|\displaystyle\sup_i m_{i,l}^{\vee}\ast f\|_p
\lesssim \displaystyle\sum_{k>0}\|(\displaystyle\sum_i|m_{i,l}^{\vee}\ast P_{i,l,k}\ast f|^2)^{\frac{1}{2}}\|_p\\+\displaystyle\sum_{k\leq 0}\|(\displaystyle\sum_i|(m_{i,l}^{\vee}-\psi_{i,l}^{\vee})\ast P_{i,l,k}\ast f|^2)^{\frac{1}{2}}\|_p +\|\displaystyle\sup_i\psi_{i,l}^{\vee}\ast f\|_p.
\end{multline*}

Note that there exist $(m,n)$ such that $\partial_{t_1}^m\partial_{t_2}^nP(t_1,t_2)$ is a monomial. Now, we determine that $\tilde{P}(t_1,t_2)=P(t_1,t_2)-a_{mn}t_1^mt_2^n$. Therefore, 
\begin{equation*}
    |m_{i,l}(\xi)|\lesssim \frac{1}{|2^{-(m,n)\cdot(\vec{v}i+\vec{w}l)}\xi_3|^{\frac{1}{m+n}}},\ \ |m_{i,l}(\xi)-\psi_{i,l}(\xi)|\lesssim|2^{(m,n)\cdot(\vec{v}i+\vec{w}l)}\xi_3|
\end{equation*}
by the van der Corput lemma and mean value theorem. Therefore, the desired result is reduced to showing that \begin{equation*}
    \|\displaystyle\sup_i\psi_{i,l}^{\vee}\ast f\|_p< C\|f\|_p\ \textrm{for}\ p>1\ \textrm{and C does not depend on}\ j.
\end{equation*}
Indeed, we deal with the remaining part using classical arguments, such bootstrap arguments and the interpolation theorem of vector-valued inequalities. As a result, the desired result follows from induction on the number of monomials of $P(t_1,t_2)$. We omit the details of the proof.

$$\ $$

Let us temporally denote $ \delta_{\vec{v}_k}(t_1^{m_s}), \delta_{\vec{v}_k}(t_1^{M_k}t_2^{N_k}), \delta^{(M_k,N_k)}  $($m_s,M_k,N_k\geq 2$) by $\delta_1,\delta_2,\delta$, for notational simplicity.

\begin{lem} For a given $\vec{v}_k$ in $\mathbb{Z}_+^2$, 
$\alpha(p)=l^{\vec{v_k}}(\frac{1}{p})=(1-\frac{\delta_1}{\delta_2})\frac{1}{p}+\frac{1}{\delta_2}$ passes through

$$(\frac{1}{p},\alpha(p))= (\frac{1}{\delta}-\frac{1}{N_k},\frac{1}{\delta}),(\frac{1}{\delta_1+\delta_2},\frac{2}{\delta_1+\delta_2}),(\frac{1}{\delta_1},\frac{1}{\delta_1}).$$
\end{lem}

$Proof$  Note that $$\frac{1}{\delta}-\frac{1}{N_k}=
\frac{1}{\delta}(\frac{N_k-\delta}{N_k})=\frac{1}{\delta}(\frac{\delta_2-\delta}{\delta_2-\delta_1}),$$ by similarity (see Figure 3).
Thus, $$(1-\frac{\delta_1}{\delta_2})(\frac{1}{\delta}-\frac{1}{N_k})+\frac{1}{\delta_2}=(1-\frac{\delta_1}{\delta_2})(\frac{1}{\delta}(\frac{\delta_2-\delta}{\delta_2-\delta_1}))+\frac{1}{\delta_2}=\frac{1}{\delta}.$$
Furthermore, we can verify directly that $\alpha(p)=(1-\frac{\delta_1}{\delta_2})\frac{1}{p}+\frac{1}{\delta_2}$ passes through
$ (\frac{1}{\delta_1+\delta_2},\frac{2}{\delta_1+\delta_2}), (\frac{1}{\delta_1},\frac{1}{\delta_1})$.
\begin{figure*}[ht]
  \includegraphics[width=0.65\textwidth]{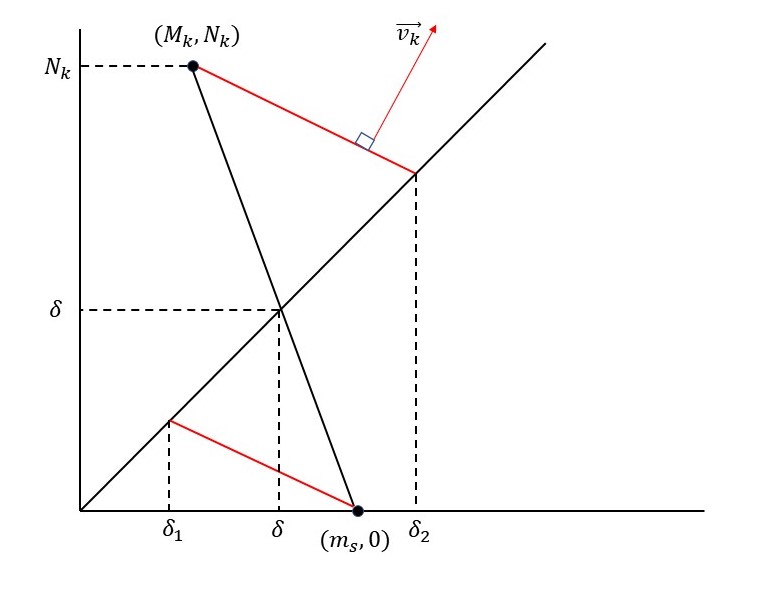}
\caption{Similarity}  
\label{3}
\end{figure*}

\begin{rmk}
By the same argument, $\alpha(p)=l^{\vec{v_{k-1}}}( \frac{1}{p})$ also passes through $(\frac{1}{\delta}-\frac{1}{N_k},\frac{1}{\delta})$. As a result, by comparing the slopes of $l^{\vec{v}_{k}}(\frac{1}{p})$ and $l^{\vec{v}_{k-1}}(\frac{1}{p})$, we can observe that
\begin{align*}
l^{\vec{v_{k-1}}}(\frac{1}{p})&>l^{\vec{v_{k}}}(\frac{1}{p}),\ \textrm{where}\  \frac{1}{p}<\frac{1}{\delta}-\frac{1}{N_k}\\
l^{\vec{v_{k-1}}}(\frac{1}{p})&<l^{\vec{v_{k}}}(\frac{1}{p}),\ \textrm{where}\  \frac{1}{p}>\frac{1}{\delta}-\frac{1}{N_k}
\end{align*}
\end{rmk}
$$\ $$
Let us temporally denote $\delta_{\vec{v}}(t_1^m),\delta_{\vec{v}}(t_1^Mt_2^N),\delta_{\vec{w}}(t_1^m),\delta_{\vec{w}}(t_1^Mt_2^N)$ ($m,M,N\geq 2$) by $\delta_{\vec{v}}^1,\delta_{v}^2,\delta_{\vec{w}}^1,\delta_{\vec{w}}^2$, for notational simplicity.

\begin{lem} Let $\vec{v}$ and $\vec{w}$ be vectors in $\mathbb{Z}_+^2$. Assume that the slope of $\vec{w}$ is greater than that of $\vec{v}$, and $\delta_{\vec{v}}^1\leq\delta_{v}^2$. Then, 
\begin{align*}
&(1)\ \delta_{\vec{v}}^2\delta_{\vec{w}}^1-\delta_{\vec{v}}^1\delta_{\vec{w}}^2<0\\
&(2)\ (M,N)\in\mathcal{R}\ \textrm{if and only if}\\ &\qquad\frac{1}{\delta_{\vec{w}}^1+\delta_{\vec{w}}^2}\leq\frac{\delta_{\vec{w}}^2-\delta_{\vec{v}}^2}{\delta_{\vec{v}}^1\delta_{\vec{w}}^2-\delta_{\vec{v}}^2\delta_{\vec{w}}^1}
\textrm{and}\ \frac{\delta_{\vec{w}}^2-\delta_{\vec{v}}^2}{\delta_{\vec{v}}^1\delta_{\vec{w}}^2-\delta_{\vec{v}}^2\delta_{\vec{w}}^1}=\frac{1}{\delta^{(M,N)}}-\frac{1}{N}
\end{align*}

\end{lem}

$\ $

$Proof$ (1) In the case that $M\leq N$,
as the slope of $\vec{w}$ is greater than that of $\vec{v}$, $\delta_{\vec{v}}^2\leq\delta_{\vec{w}}^2$ and $\delta_{\vec{w}}^1<\delta_{\vec{v}}^1$. Thus, $\delta_{\vec{v}}^2\delta_{\vec{w}}^1-\delta_{\vec{v}}^1\delta_{\vec{w}}^2\leq \delta_{\vec{w}}^2\delta_{\vec{w}}^1-\delta_{\vec{v}}^1\delta_{\vec{w}}^2=\delta_{\vec{w}}^2(\delta_{\vec{w}}^1-\delta_{\vec{v}}^1)<0$.
In the case that $M>N$, by similarity (see above Figure 4), $\frac{m}{M'}=\frac{\delta_{\vec{w}}^1}{\delta}=\frac{\delta_{\vec{v}}^1}{\delta_{\vec{v}}^2}.$
Thus, $\delta_{\vec{v}}^2\delta_{\vec{w}}^1-\delta_{\vec{v}}^1\delta_{\vec{w}}^2<\delta_{\vec{v}}^2\delta_{\vec{w}}^1-\delta_{\vec{v}}^1\delta=0$, as $\delta_{\vec{w}}^2>\delta$.
$\ $

$\ $

(2) By a simple computation, it holds that $\frac{1}{\delta_{\vec{w}}^1+\delta_{\vec{w}}^2}\leq\frac{\delta_{\vec{w}}^2-\delta_{\vec{v}}^2}{\delta_{\vec{v}}^1\delta_{\vec{w}}^2-\delta_{\vec{w}}^2\delta_{\vec{w}}^1}$ if and only if $\delta_{\vec{v}}^1-\delta_{\vec{w}}^1 \leq\delta_{\vec{w}}^2-\delta_{\vec{v}}^2$. Since $\delta_{\vec{w}}^2<\delta_{v}^2$ and $\delta_{\vec{w}}^1<\delta_{v}^1$ when $M> N$, it is necessary for $\delta_{\vec{v}}^1-\delta_{\vec{w}}^1 \leq\delta_{\vec{w}}^2-\delta_{\vec{v}}^2$ that $M\leq N$. Thus,
$$\delta_{\vec{v}}^1-\delta_{\vec{w}}^1 \leq\delta_{\vec{w}}^2-\delta_{\vec{v}}^2\Leftrightarrow \delta^{(M,N)}\leq N-\delta^{(M,N)}\ \textrm{(See Figure 4)},$$ which implies that $(M,N)\in\mathcal{R}$.

Finally, by similarity (see Figure 4),
\begin{figure*}[h]
  \includegraphics[width=0.49\textwidth]{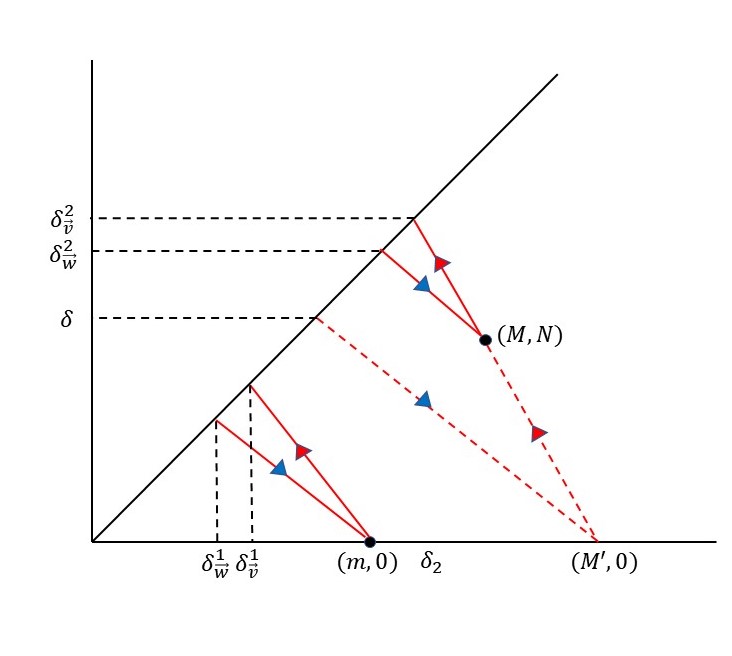}
  \label{4}
   \includegraphics[width=0.5\textwidth]{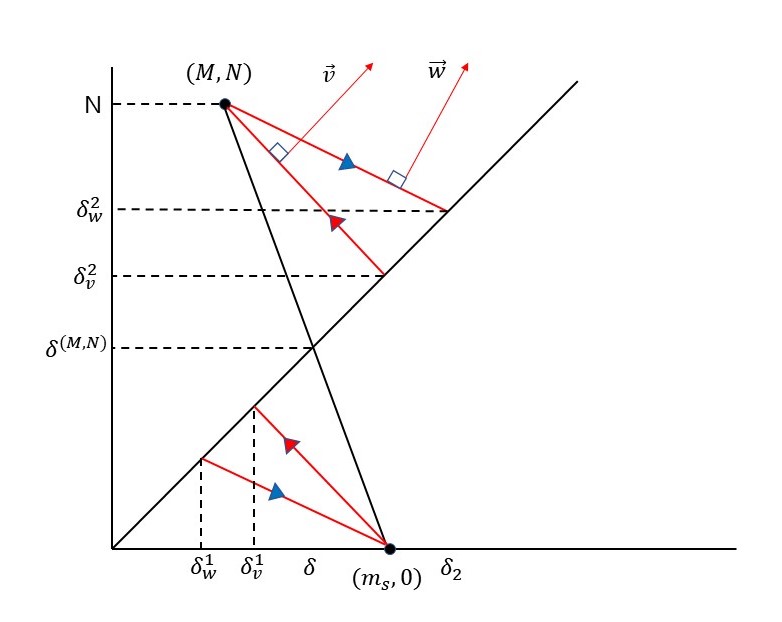}
\label{10}
\caption{Similarity}  
\end{figure*}
$$\frac{\delta_{w}^2-\delta_{v}^2}{\delta_{v}^1\delta_{w}^2-\delta_{v}^2\delta_{w}^1}=\frac{\delta_{w}^2-\delta_{v}^2}{\delta_{v}^1(\delta_{w}^2-\delta_{v}^2)+\delta_{v}^2(\delta_v^1-\delta_{w}^1)}=\frac{N-\delta^{(M,N)}}{\delta_{v}^1(N-\delta^{(M,N)})+\delta_{v}^2\delta^{(M,N)}}$$
and $\frac{\delta^{(M,N)}-\delta_v^1}{\delta_v^2-\delta_v^1}=\frac{\delta^{(M,N)}}{N}(\Leftrightarrow\delta_{v}^1(N-\delta^{(M,N)})+\delta_{v}^2\delta^{(M,N)}=\delta^{(M,N)}N)$. Therefore, $\frac{\delta_{w}^2-\delta_{v}^2}{\delta_{v}^1\delta_{w}^2-\delta_{v}^2\delta_{w}^1}=\frac{1}{\delta^{(M,N)}}-\frac{1}{N}$. In the same manner, we can observe that $\frac{1}{\delta_{\vec{w}}^1+\delta_{\vec{w}}^2}\leq\frac{\delta_{\vec{w}}^2-\delta_{\vec{v}}^2}{\delta_{\vec{v}}^1\delta_{\vec{w}}^2-\delta_{\vec{v}}^2\delta_{\vec{w}}^1} $ can be replaced by $\frac{1}{\delta_{\vec{v}}^1+\delta_{\vec{v}}^2}\leq\frac{\delta_{\vec{w}}^2-\delta_{\vec{v}}^2}{\delta_{\vec{v}}^1\delta_{\vec{w}}^2-\delta_{\vec{v}}^2\delta_{\vec{w}}^1} $ in (2).

Let us temporally denote $\delta_{\vec{v}}(t_1^m),\delta_{\vec{v}}(t_1^Mt_2^N)$  ($m,M,N\geq 2$) by $\delta_{\vec{v}}^1,\delta_{v}^2$, for notational simplicity.
\begin{prop}
Let $\vec{v}$ be a vector in $\mathbb{Z}_+^2$, and define
$$\mathcal{A}_I(f)(x)=\int f(x_1-t_1,x_2-t_2,x_3-P(t_1,t_2))\displaystyle\sum_{(j_1,j_2)\in I}\eta(\frac{t_1}{2^{-j_1}},\frac{t_2}{2^{-j_2}})dt_1dt_2,$$
where $\Lambda(P)=\{(m,0),(M,N)\} (m,M,N\geq2)$ and $I= \{\vec{v}i|i\in \mathbb{Z}_+\}$. Assume that $\delta_{\vec{v}}(t_1^m)\neq\delta_{\vec{v}}(t_1^Mt_2^N)$. Then, the operator $\mathcal{A}_I$ maps $L^p$ to $L^p_{\alpha}$ for $\alpha<\alpha(p)$ under the following conditions:

$\ $
\begin{align*}
&\bullet \textrm{In the case that}\ \delta_1< \delta_2  \\
&\alpha(p):=\left\{\begin{array}{l}        \frac{1}{\delta_1}\quad\qquad\qquad\quad\ \  \textrm{if}\  min(\frac{1}{\delta_1},\frac{1}{2})<\frac{1}{p}<\frac{1}{2}\\ (1-\frac{\delta_1}{\delta_2})\frac{1}{p}+\frac{1}{\delta_2}\ \quad\textrm{if}\  \frac{1}{\delta_1+\delta_2}<\frac{1}{p}<min(\frac{1}{\delta_1},\frac{1}{2})\\
\frac{2}{p} \quad\qquad\qquad\qquad \textrm{if}\ \frac{1}{p} \leq \frac{1}{\delta_1+\delta_2}\end{array}\right. \\
   &\bullet \textrm{In the case that}\ \delta_1>\delta_2\\
  &\alpha(p):=\left\{\begin{array}{l} \frac{1}{\delta_2}\quad\qquad\textrm{if}\  \frac{1}{2\delta_2}<\frac{1}{p}<\frac{1}{2}\\
\frac{2}{p} \quad\qquad\ \textrm{if}\ \frac{1}{p} \leq \frac{1}{2\delta_2}\end{array}\right.
\end{align*}

\end{prop}

$Proof$ In the case that $\delta_1<\delta_2$, we may assume that $\frac{1}{\delta_1}<\frac{1}{2}$, because the case with $\frac{1}{\delta_1}>\frac{1}{2}$ follows from the same argument.

 $\ $

It suffices to show that
\begin{align}
\|m_I^{\vee}\ast f\|_{L^p_{\alpha}}&\lesssim \|f\|_{p} \ \textrm{for} \ \alpha=\frac{1}{\delta_1},  \ \ \textrm{if}\ \ \frac{1}{\delta_1}<\frac{1}{p}\leq \frac{1}{2} \\
\|m_I^{\vee}\ast f\|_{L^p_{\alpha}}&\lesssim\|f\|_p\ \textrm{for}\ \alpha < (1-\frac{\delta_1}{\delta_2})\frac{1}{p}+\frac{1}{\delta_2},\ \textrm{if}\  \frac{1}{\delta_1+\delta_2}<\frac{1}{p}<\frac{1}{\delta_1}\\
\|m_I^{\vee}\ast f\|_{L^p_{\alpha}}&\lesssim \|f\|_{p} \ \textrm{for} \ \alpha=\frac{1}{\delta_1},  \ \ \textrm{if}\ \ \frac{1}{p}<\frac{1}{\delta_1+\delta_2}\ 
\end{align}
$\ $

For (3.1), let us define a Littlewood--Paley projection
$$\widehat{P_{i,k}\ast f}(\xi)=\chi^2(2^{-(\vec{v}\cdot(m,0)i+k)}\xi_3)=\chi^2(2^{-(\delta_1\vec{\mathbf{1}}\cdot\vec{v}i+k)}\xi_3)$$
$$(\tilde{P}_{i,k}\ast f)^{\wedge}(\xi)=|2^{-(\delta_1\vec{\mathbf{1}}\cdot\vec{v}i+k)}\xi_3|^{\frac{1}{\delta_1}}\chi^2(2^{-(\delta_1\vec{\mathbf{1}}\cdot\vec{v}i+k)}\xi_3),$$
 where $\chi$ is supported on $[-1,-\frac{1}{2}]\cup[\frac{1}{2},1]$ and $\displaystyle\sum_{k}\chi^2(2^{-k}\xi_3)=1$.
 Then,
\begin{align*}
\|&\displaystyle\sum_{i}2^{-\vec{v}\cdot\vec{\mathbf{1}}i}m_{i}^{\vee}\ast f\|_{L^p_{\alpha}}\\
&\lesssim\displaystyle\displaystyle\sum_{k>0}\|\displaystyle\displaystyle\sum_{i}2^{-\vec{v}\cdot\vec{\mathbf{1}}i}m_{i}^{\vee}\ast P_{i,k}\ast f\|_{L^p_{\alpha}}+\displaystyle\sum_{k\leq0}\|\displaystyle\displaystyle\sum_{i}2^{-\vec{v}\cdot\vec{\mathbf{1}}i}m_{i}^{\vee}\ast P_{i,k}\ast f\|_{L^p_{\alpha}}\\
&\lesssim\displaystyle\displaystyle\sum_{k>0}2^{\frac{1}{\delta_1}k}\|\displaystyle\displaystyle\sum_{i}m_{i}^{\vee}\ast \tilde{P}_{i,k}\ast f\|_{p}+\displaystyle\sum_{k\leq0}2^{\frac{1}{\delta_1}k}\|\displaystyle\displaystyle\sum_{i}m_{i}^{\vee}\ast \tilde{P}_{i,k}\ast f\|_{p}.
\end{align*}

By Lemma 1 and the van der Corput lemma for the variable $t_1$, we have that
$m_i(\xi)=O(|\frac{\xi_3}{2^{\vec{v}\cdot(m,0)i}}|^{-\frac{1}{2}})=O(|\frac{\xi_3}{2^{\delta_1\vec{\mathbf{1}}\cdot\vec{v}i}}|^{-\frac{1}{2}})$ for large $i$. Thus, 
\begin{align}
&\|\displaystyle\displaystyle\sum_{i}m_{i}^{\vee}\ast \tilde{P}_{i,k}\ast f\|_{L^2}\lesssim 2^{-\frac{1}{2}k}\|f\|_{L^2}\\
&\|\displaystyle\displaystyle\sum_{i}m_{i}^{\vee}\ast \tilde{P}_{i,k}\ast f\|_{p_0}\lesssim \|f\|_{p_0}\ \textrm{for}\ p_0>1,
\end{align}
$\ $

by the Littlewood--Paley theorem and Lemma 2.

$\ $

Then, by interpolation between (3.4) and (3.5), we obtain
$$\|\displaystyle\displaystyle\sum_{i}m_{i}^{\vee}\ast \tilde{P}_{i,k}\ast f\|_{p}\lesssim 2^{-(\frac{1}{p}-\epsilon)k}\|f\|_{p}\ \textrm{for}\ p>1$$.

Therefore, 
$\displaystyle\displaystyle\sum_{k>0}2^{\frac{1}{\delta_1}k}\|\displaystyle\displaystyle\sum_{i}m_{i}^{\vee}\ast \tilde{P}_{i,k}\ast f\|_{p}\lesssim \|f\|_p\ \textrm{for}\ \frac{1}{\delta_1}<\frac{1}{p}\leq \frac{1}{2}$.

Similarly, by using $m_I(\xi)=O(1)$ instead of $O(|\frac{\xi_3}{2^{\delta_1\vec{\mathbf{1}}\cdot\vec{v}i}}|^{-\frac{1}{2}})$, we have that $\displaystyle\sum_{k\leq0}2^{\frac{1}{\delta_1}k}\|\displaystyle\displaystyle\sum_{i}m_{i}^{\vee}\ast \tilde{P}_{i,k}\ast f\|_{p}\lesssim \|f\|_p\ \textrm{for}\ p>1$, which implies (3.1).

$\ $

$\ $

To prove (3.2), let us define a Littlewood--Paley projection

$$\widehat{P_{i,k}\ast f(\xi)}=\chi^2(2^{-(\vec{v}\cdot(M,N)i+k)}\xi_3)\hat{f}(\xi)=\chi^2(2^{-(\delta_2\vec{\mathbf{1}}\cdot{\vec{v}}i+k)}\xi_3)\hat{f}(\xi)$$
$$(\tilde{P_{i,k}}\ast f)^{\wedge}(\xi)=|2^{-(\delta_2\vec{\mathbf{1}}\cdot\vec{v}i+k)}\xi_3|^{\frac{1}{\delta_2}}\chi^2(2^{-(\delta_2\vec{\mathbf{1}}\cdot\vec{v}i+k)}\xi_3)\hat{f}(\xi).$$
$\ $
We estimate $\|m_I^{\vee}\ast f\|_{L^p_{\alpha}}$ as follows:
\begin{align*}
\|&\displaystyle\sum_i2^{-\vec{v}\cdot\vec{1}i}m_i^{\vee}\ast f\|_{L^p_{\alpha}}\\&\lesssim \displaystyle\sum_{k=-\infty}^{0}\|\displaystyle\displaystyle\sum_i2^{-\vec{v}\cdot\vec{\mathbf{1}}i}m_i^{\vee}\ast P_{i,k}\ast f\|_{L^p_{\alpha}} + \displaystyle\sum_{k=0}^{\infty}\|\displaystyle\displaystyle\sum_i2^{-\vec{v}\cdot\vec{\mathbf{1}}i}m_i^{\vee}\ast P_{i,k}\ast f\|_{L^p_{\alpha}}\\
&\lesssim \displaystyle\sum_{k=-\infty}^{0}2^{\frac{1}{\delta_2}k}\|\displaystyle\displaystyle\sum_im_i^{\vee}\ast \tilde{P}_{i,k}\ast f\|_{L^p_{\beta}}+ \displaystyle\sum_{k=0}^{\infty}2^{\frac{1}{\delta_2}k}\|\displaystyle\displaystyle\sum_im_i^{\vee}\ast \tilde{P}_{i,k}\ast f\|_{L^p_{\beta}},
\end{align*}
$\ $

 where $\beta<\frac{1}{p}(1-\frac{\delta_1}{\delta_2})$.

$\ $

Recall that $m_i(\xi)=\displaystyle\int e^{i\xi_3\phi_I(t_1,t_2,\xi)}\eta(t_1,t_2)dt$. We estimate the determinant of the mixed Hessian of $\phi_I(t_1,t_2,\xi)$. Owing to Lemma 1, we have that
$$ |\textrm{det}(\nabla_{t_1t_2}^2\phi_I(t_1,t_2,\xi))|\sim|2^{-\vec{v}\cdot(m,0)i}||2^{-\vec{v}\cdot(M,N)i}|\sim|2^{-\delta_1\vec{\mathbf{1}}\cdot\vec{v}i}||2^{-\delta_2\vec{\mathbf{1}}\cdot\vec{v}i}|$$ for large $i$. 
Thus, $m_i(\xi)$ =$O(|\frac{\xi_3}{2^{\delta_1\vec{\mathbf{1}}\cdot\vec{v}i}}|^{-\frac{1}{2}}|\frac{\xi_3}{2^{\delta_2\vec{\mathbf{1}}\cdot\vec{v}i}}|^{-\frac{1}{2}})$.

Furthermore, note that
$$|\frac{\xi_3}{2^{\delta_1\vec{\mathbf{1}}\cdot\vec{v}i}}|^{\frac{1}{2}}=|\frac{\xi_3}{2^{\delta_2\vec{\mathbf{1}}\cdot\vec{v}i}}|^{\frac{\delta_1}{2\delta_2}}|\xi_3|^{\frac{1}{2}-\frac{\delta_1}{2\delta_2}}.$$

Thus,
\begin{align}
\|\displaystyle\displaystyle\sum_im_i^{\vee}\ast\tilde{P}_{i,k}\ast f\|_{L^2_s}&\lesssim 2^{(1-s)k}\|f\|_{L^2}
,\textrm{where}\  s=\frac{1}{2}(1-\frac{\delta_1}{\delta_2})\\
\|\displaystyle\displaystyle\sum_im_i^{\vee}\ast\tilde{P}_{i,k}\ast f\|_{p_0}&\lesssim \|f\|_{p_0}\ \textrm{for}\ p_0>1,
\end{align}
by Lemma 2 and the Littlewood--Paley theorem.
Then, by interpolation between (3.6) and (3.7), we obtain
$$\ \|\displaystyle\displaystyle\sum_im_i^{\vee}\ast \tilde{P}_{i,k}\ast f\|_{L^p_{\beta}}\lesssim 2^{-(\frac{2}{p}-\epsilon)(1-s)k}\|f\|_{p}\ \ \textrm{for}\ 2<p<\infty,\ \textrm{where}\ \beta<\frac{2}{p}s.$$
Therefore,
\begin{equation*}
\displaystyle\sum_{k=0}^{\infty}2^{\frac{1}{\delta_2}k}\|\displaystyle\displaystyle\sum_im_i^{\vee}\ast \tilde{P}_{i,k}\ast f\|_{L^p_{\beta}}
\lesssim\displaystyle\sum_{k=0}^{\infty}2^{\frac{1}{\delta_2}k} 2^{-(\frac{2}{p}-\epsilon)(1-s)k}\|f\|_{L^p}\lesssim\|f\|_{L^p}\ \ \ \textrm{for}\ \frac{1}{p}>\frac{1}{\delta_1+\delta_2}.
\end{equation*}
Similarly, by using $m_i(\xi)=O(|\frac{1}{2^{\delta_1\vec{\mathbf{1}}\cdot\vec{v}i}}|^{-\frac{1}{2}})$ instead of $O(|\frac{\xi_3}{2^{\delta_1\vec{\mathbf{1}}\cdot\vec{v}i}}|^{-\frac{1}{2}}|\frac{\xi_3}{2^{\delta_2\vec{\mathbf{1}}\cdot\vec{v}i}}|^{-\frac{1}{2}})$, we obtain

\begin{equation*}
\displaystyle\sum_{k=-\infty}^{0}2^{\frac{1}{\delta_2}k}\|\displaystyle\displaystyle\sum_im_I^{\vee}\ast \tilde{P}_{i,k}\ast f\|_{L^p_{\beta}}
\lesssim\displaystyle\sum_{k=-\infty}^{0}2^{\frac{1}{\delta_2}k} 2^{-(\frac{2}{p}-\epsilon)(\frac{1}{2}-s)k}\|f\|_{L^p}\lesssim\|f\|_{L^p}\ \textrm{for}\ \frac{1}{p}<\frac{1}{\delta_1}
\end{equation*}

which implies (3.2).
Furthermore, interpolation between the above results and the fact that $\|\mathcal{A}_I\|_{BMO}\lesssim \|f\|_{\infty}$ yields the result for the range $\frac{1}{p}\leq\frac{1}{\delta_1+\delta_2}$.

$\ $

In the case that $\delta_1>\delta_2$, 
note that 

\begin{equation*}
\textrm{det} \begin{bmatrix}
\partial_{t_1}^2(t_1^Mt_2^N)&\partial_{t_1t_2}(t_1^Mt_2^N)\\
\partial_{t_2t_1}(t_1^Mt_2^N)&\partial_{t_2}^2(t_1^Mt_2^N)
\end{bmatrix}=MN(1-M-N)t_1^{2(M-1)}t_2^{2(N-1)}\neq 0
\ 
\end{equation*}

$\textrm{for}\  |(t_1,t_2)|\neq 0.$

Then, by Lemma 1, we have that
$m_i(\xi)=O(|\frac{\xi_3}{2^{\delta_2\vec{\mathbf{1}}\cdot\vec{v}i}}|^{-1})\ \ \textrm{for large}\ i$. Thus, it follows from the same argument used in the proof of (3.1) that

$\|m_i^{\vee}\ast f\|_{L^p_{\alpha}}\lesssim \|f\|_{p} \ \textrm{for} \ \alpha=\frac{1}{\delta_2}  \ \ \textrm{if}\ \ \frac{1}{2\delta_2}<\frac{1}{p}\leq \frac{1}{2}$. Furthermore, interpolation between the above results and the fact that $\|\mathcal{A}_I\|_{BMO}\lesssim \|f\|_{\infty}$ yields the result for the range $\frac{1}{p}\leq\frac{1}{2\delta_2}$.

$\ $

\begin{rmk}

Let $\vec{v}$ be a vector in $\mathbb{Z}_+^2$. Let $ P(t_1,t_2)$ correspond to $$\Lambda(P)=\{(m_1,0),(m_2,0),..,(m_l,0),(M_1,N_1),..,(M_n,N_n)\}.$$ Furthermore, let $(M_i,N_i)$ and $m_i$ be enumerated such that $N_1\geq N_2\geq\cdots\geq N_n$ and $m_1< m_2<\cdots< m_l$. Assume that 
$\overline{(M_i,N_i)(M_{i+1},N_{i+1})}$ is orthogonal to $\vec{v}$ for all $i$ and $m_i, M_i,N_i\geq2$.

$\ $

$(1)$ $\delta_{\vec{v}}(t_1^{m_1})<\delta_{\vec{v}}(t_1^{M_1}t_2^{N_1})$

Assume that $\partial^2_{t_2}(\displaystyle\sum_ia_{M_iN_i}t_1^{M_i}t_2^{N_i})\neq0 $ for $|(t_1,t_2)|\neq0$. Then,
 $\mathcal{A}_I$ for $P(t_1,t_2)$ yields the same result as in Proposition 1. Indeed, we obtain $$|\textrm{det}(\nabla_{t_1t_2}^2\phi_I(t_1,t_2,\xi))|\sim|2^{-\delta_1\vec{\mathbf{1}}\cdot\vec{v}i}||2^{-\delta_2\vec{\mathbf{1}}\cdot\vec{v}i}|\ \textrm{for large}\ i, $$  where $\delta_1=\delta_{\vec{v}}(t_1^{m_1})$ and $\delta_2=\delta_{\vec{v}}(t_1^{M_1}t_2^{N_1})$. 
 This implies that $$m_I(\xi)=O(|\frac{\xi_3}{2^{\delta_1\vec{\mathbf{1}}\cdot\vec{v}i}}|^{-\frac{1}{2}}|\frac{\xi_3}{2^{\delta_2\vec{\mathbf{1}}\cdot\vec{v}i}}|^{-\frac{1}{2}}),$$ as in Proposition 1. 
 
 Note that even if the `real' second dominant monomial in view of Lemma 1 is $t_1^{m_2}$, this monomial cannot affect the value of the determinant of he mixed Hessian of $\phi_{I}(t_1,t_2,\xi)$, because $\partial_{t_2}t_1^{m_2}=0$. Thus, we choose our second dominant monomials among those of the form $t_1^Mt_2^N$ ($M,N \neq 0$).

$\ $

$(2)$ $\delta_{\vec{v}}(t_1^{m_1})> \delta_{\vec{v}}(t_1^{M_1}t_2^{N_1})$

Assume that $| det(\nabla^2_{t_1t_2} (\displaystyle\sum_ia_{M_iN_i}t_1^{M_i}t_2^{N_i}))|\neq 0$ for $|(t_1,t_2)|\neq0$. Then,
 $\mathcal{A}_I$ for $P(t_1,t_2)$ yields the same result as in Proposition 1. Indeed, 
we obtain $$|\textrm{det}(\nabla_{t_1t_2}^2\phi_I(t_1,t_2,\xi))|\sim|2^{-\vec{v}\cdot(M_1,N_1)i}|^2\ \textrm{for large}\ i,$$ which implies that $m_I(\xi)=O(|\frac{\xi_3}{2^{\vec{v}\cdot(M_1,N_1)i}}|^{-1})$, as in Proposition 1.

$\ $

$(3)$ $\delta_{\vec{v}}(t_1^{m_1})=\delta_{\vec{v}}(t_1^{M_1}t_2^{N_1})$

 Assume that $| det(\nabla^2_{t_1t_2}(t_1^{m_1}+\displaystyle\sum_ia_{M_iN_i}t_1^{M_i}t_2^{N_i}))|\neq 0$ for $|(t_1,t_2)|\neq0$. Then, $\mathcal{A}_I$ for $P(t_1,t_2)$ yields the same result as in Proposition 1 for the case that $\delta_1>\delta_2$. Indeed, we obtain 
 $$|\textrm{det}(\nabla_{t_1t_2}^2\phi_I(t_1,t_2,\xi))|\sim |2^{-\vec{v}\cdot(M_1,N_1)i}|^2\ \textrm{for large}\ i,$$ which implies that $m_I(\xi)=O(|\frac{\xi_3}{2^{\vec{v}\cdot(M_1,N_1)i}}|^{-1})$, as in Proposition 1.
\end{rmk}

$$\ $$

Let us temporally denote $\delta_{\vec{v}}(t_1^m),\delta_{\vec{v}}(t_1^Mt_2^N),\delta_{\vec{w}}(t_1^m),\delta_{\vec{w}}(t_1^Mt_2^N)$ by $\delta_{\vec{v}}^1,\delta_{v}^2,\delta_{\vec{w}}^1,$ $\delta_{\vec{w}}^2$, for notational simplicity.

\begin{prop}

Let $\vec{v}$ and $\vec{w}$ be vectors in $\mathbb{Z}_+^2$, 
and define
$$\mathcal{A}_J(f)(x)=\int f(x_1-t_1,x_2-t_2,x_3-P(t_1,t_2))\displaystyle\sum_{(j_1,j_2)\in J}\eta(\frac{t_1}{2^{-j_1}},\frac{t_2}{2^{-j_2}})dt_1dt_2,$$
 where $\Lambda(P)=\{(m,0),(M,N)\} (m,M,N \geq 2)$ and $J= \{\vec{v}i+\vec{w}l|i,l\in \mathbb{Z}_+\}$. 
Assume that the slope of $\vec{w}$ is greater than that of $\vec{v}$, and $\vec{v}$ and $\vec{w}$ satisfy either

$$\delta_{\vec{v}}^1<\delta_{\vec{v}}^2,\ \delta_{\vec{w}}^1<\delta_{\vec{w}}^2\ \textrm{or}\ \delta_{\vec{v}}^1>\delta_{\vec{v}}^2,\ \ \delta_{\vec{w}}^1>\delta_{\vec{w}}^2.$$ Then, the operator $\mathcal{A}_J$ maps $L^p$ to $L^p_{\alpha}$ for $\alpha<\alpha(p)$ under the following conditions. 

$\ $

$\bullet \textrm{In the case that}\ (M,N)\not\in\mathcal{R}$

\begin{align*}
 1)\ & \delta_{\vec{v}}^1<\delta_{\vec{v}}^2,\ \delta_{\vec{w}}^1<\delta_{\vec{w}}^2   \\
\alpha(p):=&\left\{\begin{array}{l} \frac{1}{\delta_{\vec{v}}^1}\quad\qquad\textrm{if}\  min(\frac{1}{\delta_{\vec{v}}^1},\frac{1}{2})<\frac{1}{p}<\frac{1}{2}\\ (1-\frac{\delta_{\vec{v}}^1}{\delta_{\vec{v}}^2})\frac{1}{p}+\frac{1}{\delta_{\vec{v}}^2}\ \textrm{if}\  \frac{1}{\delta_{\vec{v}}^1+\delta_{\vec{v}}^2}<\frac{1}{p}<min(\frac{1}{\delta_{\vec{v}}^1},\frac{1}{2})\\
\frac{2}{p} \quad\qquad \textrm{if}\ \frac{1}{p} \leq \frac{1}{\delta_{\vec{v}}^1+\delta_{\vec{v}}^2}\end{array}\right.  \\
2)\ & \delta_{\vec{v}}^1>\delta_{\vec{v}}^2,\ \delta_{\vec{w}}^1>\delta_{\vec{w}}^2 \\
    \alpha(p):=&\left\{\begin{array}{l} min(\frac{1}{\delta_{\vec{w}}^2}, \frac{1}{\delta_{\vec{v}}^2})\qquad\textrm{if}\  min(\frac{1}{2\delta_{\vec{w}}^2}, \frac{1}{2\delta_{\vec{v}}^2})<\frac{1}{p}<\frac{1}{2}\\
\frac{2}{p} \quad\qquad \textrm{if}\ \frac{1}{p} \leq min(\frac{1}{2\delta_{\vec{w}}^2}, \frac{1}{2\delta_{\vec{v}}^2})\end{array}\right.
\end{align*}

$\bullet$ In the case that $(M,N)\in\mathcal{R}$

\begin{align*}
\qquad1)\ &\delta_{\vec{v}}^1<\delta_{\vec{v}}^2,\ \ \delta_{\vec{w}}^1<\delta_{\vec{w}}^2    \\
\alpha(p):=&\left\{\begin{array}{l} \frac{1}{\delta_{\vec{v}}^1}\quad\qquad\textrm{if}\  min(\frac{1}{\delta_{\vec{v}}^1},\frac{1}{2})<\frac{1}{p}<\frac{1}{2}\\ (1-\frac{\delta_{\vec{v}}^1}{\delta_{\vec{v}}^2})\frac{1}{p}+\frac{1}{\delta_{\vec{v}}^2}\ \textrm{if}\  \frac{1}{\delta^{(M,N)}}-\frac{1}{N}<\frac{1}{p}<min(\frac{1}{\delta_{\vec{v}}^1},\frac{1}{2})\\
(1-\frac{\delta_{\vec{w}}^1}{\delta_{\vec{w}}^2})\frac{1}{p}+\frac{1}{\delta_{\vec{w}}^2}\ \textrm{if}\  \frac{1}{\delta_{\vec{w}}^1+\delta_{\vec{w}}^2}<\frac{1}{p}<\frac{1}{\delta^{(M,N)}}-\frac{1}{N}\\
\frac{2}{p} \quad\qquad \textrm{if}\ \frac{1}{p} \leq \frac{1}{\delta_{\vec{v}}^1+\delta_{\vec{v}}^2}\end{array}\right. \\
\qquad2)\  &\delta_{\vec{v}}^1>\delta_{\vec{v}}^2,\ \ \delta_{\vec{w}}^1>\delta_{\vec{w}}^2\\
    \alpha(p):=&\left\{\begin{array}{l} min(\frac{1}{\delta_{\vec{w}}^2}, \frac{1}{\delta_{\vec{v}}^2})\qquad\textrm{if}\  min(\frac{1}{2\delta_{\vec{w}}^2}, \frac{1}{2\delta_{\vec{v}}^2})<\frac{1}{p}<\frac{1}{2}\\
\frac{2}{p} \quad\qquad \textrm{if}\ \frac{1}{p} \leq min(\frac{1}{2\delta_{\vec{w}}^2}, \frac{1}{2\delta_{\vec{v}}^2})\end{array}\right.
\end{align*}

$\ $
\end{prop}

$Proof$ We only consider $i+l\geq 2M$ for sufficiently large $M$.

$\ $
\begin{align*}
&\|\displaystyle\sum_{i,l>2M}2^{-(\vec{v}i+\vec{w}l)\cdot\vec{\mathbf{1}}}m_{i,l}^{\vee}\ast f\|_{L^p_{\alpha}}\\
&\lesssim\displaystyle\sum_{j\leq M }\|\displaystyle\sum_{l\geq M}2^{-(\vec{v}i+\vec{w}l)\cdot\vec{\mathbf{1}}}m_{i,l}^{\vee}\ast f\|_{L^p_{\alpha}} +\displaystyle\sum_{l\leq M}\|\displaystyle\sum_{j\geq M}2^{-(\vec{v}i+\vec{w}l)\cdot\vec{\mathbf{1}}}m_{i,l}^{\vee}\ast f\|_{L^p_{\alpha}}+\|\displaystyle\sum_{i,l\geq M}2^{-(\vec{v}i+\vec{w}l)\cdot\vec{\mathbf{1}}}m_{i,l}^{\vee}\ast f\|_{L^p_{\alpha}}
\end{align*}

The first two terms are dealt with as in Proposition 1, and so by the remark for the Lemma 3 we can observe that the intersection the of $L^p\rightarrow L^p_{\alpha}$ ranges of the first two terms yields the desired result. Thus, it remains only to consider the last term.

$\ $

First, consider the case that $\delta_{\vec{v}}^1<\delta_{\vec{v}}^2,\ \delta_{\vec{w}}^1<\delta_{\vec{w}}^2$. We may assume that $\frac{1}{\delta_{\vec{v}}^1}<\frac{1}{2}$, as the case with $\frac{1}{\delta_{\vec{v}}^1}\geq \frac{1}{2}$ follows from the same argument.

For the range $\frac{1}{\delta_{\vec{v}}^1}<\frac{1}{p}\leq\frac{1}{2}$, let us define a Littlewood--Paley projection
$$\widehat{P_{i,l,k}\ast f}(\xi)=\chi^2(2^{-(\delta_{\vec{v}}^1\vec{\mathbf{1}}\cdot\vec{v}i+\delta_{\vec{w}}^1\vec{\mathbf{1}}\cdot\vec{w}l+k)}\xi_3)\hat{f}(\xi)$$
$$(\tilde{P}_{i,l,k}\ast f)^{\wedge}(\xi)=|2^{-(\delta_{\vec{v}}^1\vec{\mathbf{1}}\cdot\vec{v}i+\delta_{\vec{w}}^1\vec{\mathbf{1}}\cdot\vec{w}l+k)}\xi_3|^{\frac{1}{\delta_{\vec{v}}^1}}\chi^2(2^{-(\delta_{\vec{v}}^1\vec{\mathbf{1}}\cdot\vec{v}i+\delta_{\vec{w}}^1\vec{\mathbf{1}}\cdot\vec{w}l+k)}\xi_3)\hat{f}(\xi).$$
 Then, we estimate $\|m_J(\xi)\|_{L^p_{\alpha}}$ for $\alpha=\frac{1}{\delta_{\vec{v}}^1}$ as follows:
\begin{align*}
&\|\displaystyle\sum_{i,j}2^{-(\vec{v}i+\vec{w}l)\cdot\vec{\mathbf{1}}}m_{i,l}(\xi)\ast f\|_{L^p_{\alpha}}\\
&\lesssim\displaystyle\displaystyle\sum_{k>0}\|\displaystyle\displaystyle\sum_{i}\sum_{l}2^{-\vec{v}\cdot\vec{\mathbf{1}}i-\vec{w}\cdot\vec{\mathbf{1}}l}m_{i,l}^{\vee}\ast P_{i,l,k}\ast f\|_{L^p_{\alpha}}+\displaystyle\sum_{k\leq0}\|\displaystyle\displaystyle\sum_{i}\sum_{l}2^{-\vec{v}\cdot\vec{\mathbf{1}}i-\vec{w}\cdot\vec{\mathbf{1}}l}m_{i,l}^{\vee}\ast P_{i,l,k}\ast f\|_{L^p_{\alpha}}\\
&\lesssim\displaystyle\displaystyle\sum_{k>0}\sum_{l}2^{\frac{1}{\delta_{v}^1}k}2^{-(1-\frac{\delta_{\vec{w}}^1}{\delta_{\vec{v}}^1})\vec{w}\cdot\vec{\mathbf{1}}l}\|\displaystyle\displaystyle\sum_{i}m_{i,l}^{\vee}\ast \tilde{P}_{i,l,k}\ast f\|_{p}+\displaystyle\sum_{k\leq0}\sum_l2^{\frac{1}{\delta_{v}^1}k}2^{-(1-\frac{\delta_{\vec{w}}^1}{\delta_{\vec{v}}^1})\vec{w}\cdot\vec{\mathbf{1}}l}\|\displaystyle\displaystyle\sum_{i}m_{i,l}^{\vee}\ast \tilde{P}_{i,l,k}\ast f\|_{p}.
\end{align*}

Because $1-\frac{\delta_{\vec{w}}^1}{\delta_{\vec{v}}^1}>0$, it follows that 
$$\ \|\displaystyle\sum_{i,l}2^{-(\vec{v}i+\vec{w}l)\cdot\vec{\mathbf{1}}}m_{i,l}(\xi)\ast f\|_{L^p_{\alpha}}\lesssim\|f\|_{L^2}\ \textrm{for}\ \frac{1}{\delta_{\vec{v}}^1}<\frac{1}{p}\leq\frac{1}{2},$$
by considering the fact that $m_{i,l}(\xi)=O(|\frac{\xi_3}{2^{\delta_{\vec{v}}^1\vec{1}\cdot\vec{v}i+\delta_{\vec{w}}^1\vec{1}\cdot\vec{w}l}}|^{-\frac{1}{2}})$ and employing the same argument as in Proposition 1.

For the range $\frac{1}{p}<\frac{1}{\delta_{\vec{v}}^1}$, let us define a Littlewood--Paley projection
$$\widehat{P_{i,l,k}\ast f}(\xi)=\chi^2(2^{-(\delta_{\vec{v}}^2\vec{\mathbf{1}}\cdot\vec{v}i+\delta_{\vec{w}}^2\vec{\mathbf{1}}\cdot\vec{w}l+k)}\xi_3)\hat{f}(\xi)$$
$$(\tilde{P}_{i,l,k}\ast f)^{\wedge}(\xi)=|2^{-(\delta_{\vec{v}}^2\vec{\mathbf{1}}\cdot\vec{v}i+\delta_{\vec{w}}^2\vec{\mathbf{1}}\cdot\vec{w}l+k)}\xi_3|^{\frac{1}{\delta_{\vec{v}}^2}}\chi^2(2^{-(\delta_{\vec{v}}^2\vec{\mathbf{1}}\cdot\vec{v}i+\delta_{\vec{w}}^2\vec{\mathbf{1}}\cdot\vec{w}l+k)}\xi_3)\hat{f}(\xi).$$
We estimate $\|m_{J}^{\vee}\ast f\|_{L^p_{\alpha}}$ for $\alpha<\frac{1}{p}(1-\frac{\delta_{v}^1}{\delta_{v}^2})+\frac{1}{\delta_{v}^2}$ as follows:
\begin{align*}
&\|m_{i,l}^{\vee}\ast f\|_{L^p_{\alpha}}\\
&\lesssim\displaystyle\displaystyle\sum_{k>0}\|\displaystyle\displaystyle\sum_{i}\sum_{l}2^{-\vec{v}\cdot\vec{\mathbf{1}}i-\vec{w}\cdot\vec{\mathbf{1}}l}m_{i,l}^{\vee}\ast P_{i,l,k}\ast f\|_{L^p_{\alpha}}+\displaystyle\sum_{k\leq0}\|\displaystyle\displaystyle\sum_{i}\sum_{l}2^{-\vec{v}\cdot\vec{\mathbf{1}}i-\vec{w}\cdot\vec{\mathbf{1}}l}m_{i,l}^{\vee}\ast P_{i,l,k}\ast f\|_{L^p_{\alpha}}\\
&\lesssim\displaystyle\displaystyle\sum_{k>0}\sum_{l}2^{\frac{1}{\delta_{\vec{v}}^2}k}2^{-(1-\frac{\delta_{\vec{w}}^2}{\delta_{\vec{v}}^2})\vec{w}\cdot\vec{\mathbf{1}}l}\|\displaystyle\displaystyle\sum_{i}m_{i,l}^{\vee}\ast \tilde{P}_{i,l,k}\ast f\|_{L^p_{\beta}}+\displaystyle\sum_{k\leq0}\sum_l2^{\frac{1}{\delta_{\vec{v}}^2}k}2^{-(1-\frac{\delta_{\vec{w}}^2}{\delta_{\vec{v}}^2})\vec{w}\cdot\vec{\mathbf{1}}l}\|\displaystyle\displaystyle\sum_{i}m_{i,l}^{\vee}\ast \tilde{P}_{i,l,k}\ast f\|_{L^p_{\beta}},
\end{align*}
$\ $

 where $\beta<\frac{1}{p}(1-\frac{\delta_{\vec{v}}^1}{\delta_{\vec{v}}^2})$. Note that 
$$m_{i,l}(\xi)=O(|\frac{\xi_3}{2^{\delta_{\vec{v}}^1\vec{\mathbf{1}}\cdot\vec{v}i+\delta_{\vec{w}}^1\vec{\mathbf{1}}\cdot\vec{w}l}}|^{-\frac{1}{2}}|\frac{\xi_3}{2^{\delta_{\vec{v}}^2\vec{\mathbf{1}}\cdot\vec{v}i+\delta_{\vec{w}}^2\vec{\mathbf{1}}\cdot\vec{w}l}}|^{-\frac{1}{2}})\ \textrm{and}$$
\begin{equation}
|\frac{\xi_3}{2^{\delta_{\vec{v}}^1\vec{\mathbf{1}}\cdot\vec{v}i+\delta_{\vec{w}}^1\vec{\mathbf{1}}\cdot\vec{w}l}}|^{\frac{1}{2}}=|\frac{\xi_3}{2^{\delta_{\vec{v}}^2\vec{\mathbf{1}}\cdot\vec{v}i+\delta_{\vec{w}}^2\vec{\mathbf{1}}\cdot\vec{w}l}}|^{\frac{\delta_{\vec{v}}^1}{2\delta_{\vec{v}}^2}}|\xi_3|^{\frac{1}{2}(1-\frac{\delta_{\vec{v}}^1}{\delta_{\vec{v}}^2})}2^{-\frac{1}{2\delta_{v}^2}(\delta_{v}^2\delta_{w}^1-\delta_{v}^1\delta_{w}^2)\vec{w}\cdot\vec{\mathbf{1}}l},
\end{equation}
by the determinant of the mixed Hessian of $\phi_{J}(t_1,t_2,\xi)$. Then, it follows from the same argument as in Proposition 1 that 
\begin{multline}
\displaystyle\displaystyle\sum_{k>0}\sum_{l}2^{\frac{1}{\delta_{\vec{v}}^2}k}2^{-(1-\frac{\delta_{\vec{w}}^2}{\delta_{\vec{v}}^2})\vec{w}\cdot\vec{\mathbf{1}}l}\|\displaystyle\displaystyle\sum_{i}m_{i,l}^{\vee}\ast \tilde{P}_{i,l,k}\ast f\|_{L^p_{\beta}}\\\lesssim\displaystyle\displaystyle\sum_{k>0}\sum_{l}2^{\frac{1}{\delta_{\vec{v}}^2}k-\frac{1}{p}(1+\frac{\delta_{\vec{v}}^1}{\delta_{\vec{v}}^2})k}2^{-(1-\frac{\delta_{\vec{w}}^2}{\delta_{\vec{v}}^2})\vec{w}\cdot\vec{\mathbf{1}}l+\frac{1}{p\delta_{v}^2}(\delta_{v}^2\delta_{w}^1-\delta_{v}^1\delta_{w}^2)\vec{w}\cdot\vec{\mathbf{1}}l}\|f\|_{L^p}\lesssim\|f\|_p
\end{multline}

$\ $

for 
$\frac{1}{p}(\delta_{v}^2\delta_{w}^1-\delta_{v}^1\delta_{w}^2)<\delta_{\vec{v}}^2-\delta_{\vec{w}}^2$ and $\frac{1}{\delta_{\vec{v}}^2+\delta_{\vec{v}}^1}<\frac{1}{p}$.

$\ $

Furthermore, from $m_{i,l}(\xi)=O(|\frac{\xi_3}{2^{\delta_{\vec{v}}^1\vec{\mathbf{1}}\cdot\vec{v}i+\delta_{\vec{w}}^1\vec{\mathbf{1}}\cdot\vec{w}l}}|^{-\frac{1}{2}})$ and (3.8), it follows that 
\begin{multline}
\displaystyle\sum_{k\leq0}\sum_l2^{\frac{1}{\delta_{\vec{v}}^2}k}2^{-(1-\frac{\delta_{\vec{w}}^2}{\delta_{\vec{v}}^2})\vec{w}\cdot\vec{\mathbf{1}}l}\|\displaystyle\displaystyle\sum_{i}m_{i,l}^{\vee}\ast \tilde{P}_{i,l,k}\ast f\|_{L^p_{\beta}}\\
\lesssim\displaystyle\displaystyle\sum_{k\leq 0}\sum_{l}2^{\frac{1}{\delta_{\vec{v}}^2}k-\frac{1}{p}\frac{\delta_{\vec{v}}^1}{\delta_{\vec{v}}^2}k}2^{-(1-\frac{\delta_{\vec{w}}^2}{\delta_{\vec{v}}^2})\vec{w}\cdot\vec{\mathbf{1}}l+\frac{1}{p\delta_{v}^2}(\delta_{v}^2\delta_{w}^1-\delta_{v}^1\delta_{w}^2)\vec{w}\cdot\vec{\mathbf{1}}l}\|f\|_{L^p}\lesssim\|f\|_p
\end{multline}

for $\frac{1}{p}(\delta_{v}^2\delta_{w}^1-\delta_{v}^1\delta_{w}^2)<\delta_{\vec{v}}^2-\delta_{\vec{w}}^2$ and $\frac{1}{p}<\frac{1}{\delta_{\vec{v}}^1}$.

$\ $

From (3.9), (3.10), and Lemma 3 part (1), it holds for $\frac{\delta_{\vec{v}}^2-\delta_{\vec{w}}^2}{\delta_{v}^2\delta_{w}^1-\delta_{v}^1\delta_{w}^2}<\frac{1}{p}$ and $\frac{1}{\delta_{\vec{v}}^2+\delta_{\vec{v}}^1}<\frac{1}{p}<\frac{1}{\delta_{\vec{v}}^1}$ that
\begin{equation*}
\|m_J^{\vee}\ast f\|_{L^p_{\alpha}}\lesssim\|f\|_p,\ \textrm{where}\ \alpha<\frac{1}{p}(1-\frac{\delta_{\vec{v}}^1}{\delta_{\vec{v}}^2})+\frac{1}{\delta_{\vec{v}}^2}.
\end{equation*}
Similarly, let us define a Littlewood--Paley projection
$$\widehat{P_{i,l,k}\ast f}(\xi)=\chi^2(2^{-(\delta_{\vec{v}}^2\vec{\mathbf{1}}\cdot\vec{v}i+\delta_{\vec{w}}^2\vec{\mathbf{1}}\cdot\vec{w}l+k)}\xi_3)\hat{f}(\xi)$$
$$(\tilde{P}_{i,l,k}\ast f)^{\wedge}(\xi)=|2^{-(\delta_{\vec{v}}^2\vec{\mathbf{1}}\cdot\vec{v}i+\delta_{\vec{w}}^2\vec{\mathbf{1}}\cdot\vec{w}l+k)}\xi_3|^{\frac{1}{\delta_{\vec{w}}^2}}\chi^2(2^{-(\delta_{\vec{v}}^2\vec{\mathbf{1}}\cdot\vec{v}i+\delta_{\vec{w}}^2\vec{\mathbf{1}}\cdot\vec{w}l+k)}\xi_3)\hat{f}(\xi).$$
We estimate $\|m_J^{\vee}\ast f\|_{L^p_{\alpha}}$ for $\alpha<\frac{1}{p}(1-\frac{\delta_{w}^1}{\delta_{w}^2})+\frac{1}{\delta_{w}^2}$ as follows:
\begin{align*}
&\|m_{i,l}^{\vee}\ast f\|_{L^p_{\alpha}}\\
&\lesssim\displaystyle\displaystyle\sum_{k>0}\|\displaystyle\displaystyle\sum_{l}\sum_{i}2^{-\vec{v}\cdot\vec{\mathbf{1}}i-\vec{w}\cdot\vec{\mathbf{1}}l}m_{i,l}^{\vee}\ast P_{i,l,k}\ast f\|_{L^p_{\alpha}}+\displaystyle\sum_{k\leq0}\|\displaystyle\displaystyle\sum_{l}\sum_{i}2^{-\vec{v}\cdot\vec{\mathbf{1}}i-\vec{w}\cdot\vec{\mathbf{1}}l}m_{i,l}^{\vee}\ast P_{i,l,k}\ast f\|_{L^p_{\alpha}}
\end{align*}
\begin{align*}
&\lesssim\displaystyle\displaystyle\sum_{k>0}\sum_{i}2^{\frac{1}{\delta_{\vec{w}}^2}k}2^{-(1-\frac{\delta_{\vec{v}}^2}{\delta_{\vec{w}}^2})\vec{v}\cdot\vec{\mathbf{1}}i}\|\displaystyle\displaystyle\sum_{l}m_{i,l}^{\vee}\ast \tilde{P}_{i,l,k}\ast f\|_{L^p_{\beta}}+\displaystyle\sum_{k\leq0}\sum_i2^{\frac{1}{\delta_{\vec{w}}^2}k}2^{-(1-\frac{\delta_{\vec{v}}^2}{\delta_{\vec{w}}^2})\vec{v}\cdot\vec{\mathbf{1}}i}\|\displaystyle\displaystyle\sum_{l}m_{i,l}^{\vee}\ast \tilde{P}_{i,l,k}\ast f\|_{L^p_{\beta}},
\end{align*}
$\ $

 where $ \beta<\frac{1}{p}(1-\frac{\delta_{\vec{w}}^1}{\delta_{\vec{w}}^2})$. Furthermore, note that
$$|\frac{\xi_3}{2^{\delta_{\vec{v}}^1\vec{\mathbf{1}}\cdot\vec{v}i+\delta_{\vec{w}}^1\vec{\mathbf{1}}\cdot\vec{w}l}}|^{\frac{1}{2}}=|\frac{\xi_3}{2^{\delta_{\vec{v}}^2\vec{\mathbf{1}}\cdot\vec{v}i+\delta_{\vec{w}}^2\vec{\mathbf{1}}\cdot\vec{w}l}}|^{\frac{\delta_{\vec{w}}^1}{2\delta_{\vec{w}}^2}}|\xi_3|^{\frac{1}{2}(1-\frac{\delta_{\vec{w}}^1}{\delta_{\vec{w}}^2})}2^{-\frac{1}{2\delta_{w}^2}(\delta_{w}^2\delta_{v}^1-\delta_{w}^1\delta_{v}^2)\vec{v}\cdot\vec{\mathbf{1}}i}$$
From the same argument as in (3.9) and (3.10), for $\frac{1}{p}<\frac{\delta_{\vec{v}}^2-\delta_{\vec{w}}^2}{\delta_{v}^2\delta_{w}^1-\delta_{v}^1\delta_{w}^2}$ and $\frac{1}{\delta_{\vec{w}}^2+\delta_{\vec{w}}^1}<\frac{1}{p}<\frac{1}{\delta_{\vec{w}}^1}$, it holds that
\begin{equation*}
\|m_J^{\vee}\ast f\|_{L^p_{\alpha}}\lesssim\|f\|_p,\ \textrm{where}\ \alpha<\frac{1}{p}(1-\frac{\delta_{\vec{w}}^1}{\delta_{\vec{w}}^2})+\frac{1}{\delta_{\vec{w}}^2}.
\end{equation*}
As a result, by considering Lemma 4 we conclude the following.

In the case that $(M,N)\not\in\mathcal{R}$ and $\frac{1}{\delta_{\vec{v}}^1}<\frac{1}{2}$, $\delta_{\vec{v}}^1>\delta_{\vec{v}}^2,$ and $\delta_{\vec{w}}^1>\delta_{\vec{w}}^2$:
\begin{align*}
&\|m_J^{\vee}\ast f\|_{L^p_{\alpha}}\lesssim\|f\|_p,\ \textrm{where}\ \ \frac{1}{\delta_{v}^1}<\frac{1}{p}<\frac{1}{2}\ \textrm{and}\ \alpha=\frac{1}{\delta_{\vec{v}}^1}\\
&\|m_J^{\vee}\ast f\|_{L^p_{\alpha}}\lesssim\|f\|_p,\ \textrm{where}\ \ \frac{1}{\delta_{v}^1+\delta_{v}^2}<\frac{1}{p}<\frac{1}{\delta_{v}^1}\ \textrm{and}\ \alpha<\frac{1}{p}(1-\frac{\delta_{\vec{v}}^1}{\delta_{\vec{v}}^2})+\frac{1}{\delta_{\vec{v}}^2}
\end{align*}

In the case that $(M,N)\in\mathcal{R}$ and $\frac{1}{\delta_{v}^1}<\frac{1}{2}$, $\delta_{\vec{v}}^1>\delta_{\vec{v}}^2,$ and $\delta_{\vec{w}}^1>\delta_{\vec{w}}^2$:
\begin{align*}
&\|m_J^{\vee}\ast f\|_{L^p_{\alpha}}\lesssim\|f\|_p,\ \textrm{where}\ \ \frac{1}{\delta_{v}^1}<\frac{1}{p}<\frac{1}{2}\ \textrm{and}\ \alpha=\frac{1}{\delta_{\vec{v}}^1} \\
&\|m_J^{\vee}\ast f\|_{L^p_{\alpha}}\lesssim\|f\|_p,\ \textrm{where}\ \ \frac{\delta_{w}^2-\delta_{v}^2}{\delta_{v}^1\delta_{w}^2-\delta_{v}^2\delta_{w}^1}<\frac{1}{p}<\frac{1}{\delta_{v}^1}\ \textrm{and}\ \alpha<\frac{1}{p}(1-\frac{\delta_{\vec{v}}^1}{\delta_{\vec{v}}^2})+\frac{1}{\delta_{\vec{v}}^2}\\
&\|m_J^{\vee}\ast f\|_{L^p_{\alpha}}\lesssim\|f\|_p,\ \textrm{where}\ \\
&\qquad\qquad\qquad\qquad\quad\ \ \frac{1}{\delta_{w}^1+\delta_{w}^2}<\frac{1}{p}<\frac{\delta_{w}^2-\delta_{v}^2}{\delta_{v}^1\delta_{w}^2-\delta_{v}^2\delta_{w}^1}\ \textrm{and}\ \alpha<\frac{1}{p}(1-\frac{\delta_{w}^1}{\delta_{w}^2})+\frac{1}{\delta_{w}^2}
\end{align*}
By Lemma 2(2), $\frac{\delta_{\vec{w}}^2-\delta_{\vec{v}}^2}{\delta_{\vec{v}}^1\delta_{\vec{w}}^2-\delta_{\vec{v}}^2\delta_{\vec{w}}^1}$ can be replaced by $\frac{1}{\delta^{(M,N)}}-\frac{1}{N}$.

Furthermore, interpolation between the above results and the fact that $\|\mathcal{A}_I\|_{BMO}\lesssim \|f\|_{\infty}$ yields the result for the range $\frac{1}{p}\leq \frac{1}{\delta_{\vec{v}}^1+\delta_{\vec{v}}^2}$. Thus, the proof is complete for the case that $\delta_{\vec{v}}^1>\delta_{\vec{v}}^2$ and $\delta_{\vec{w}}^1>\delta_{\vec{w}}^2$.

$\ $

Second, consider the case that $\delta_{\vec{v}}^1>\delta_{\vec{v}}^2$ and $\delta_{\vec{w}}^1>\delta_{\vec{w}}^2$. 
\begin{equation*}
\widehat{P_{i,l,k}\ast f}(\xi)=\chi^2(2^{-(\delta_{\vec{v}}^2\vec{\mathbf{1}}\cdot\vec{v}i+\delta_{\vec{w}}^2\vec{\mathbf{1}}\cdot\vec{w}l+k)}\xi_3)\hat{f}(\xi)
\end{equation*}
\begin{equation*}
(\tilde{P}_{i,l,k}\ast f)^{\wedge}(\xi)=|2^{-(\delta_{\vec{v}}^2\vec{\mathbf{1}}\cdot\vec{v}i+\delta_{\vec{w}}^2\vec{\mathbf{1}}\cdot\vec{w}l+k)}\xi_3|^{\alpha}\chi^2(2^{-(\delta_{\vec{v}}^2\vec{\mathbf{1}}\cdot\vec{v}i+\delta_{\vec{w}}^2\vec{\mathbf{1}}\cdot\vec{w}l+k)}\xi_3)\hat{f}(\xi).
\end{equation*}

We estimate $\|m_J(\xi)\|_{L^p_{\alpha}}$ for $\alpha<min(\frac{1}{\delta_{\vec{w}}^2},\frac{1}{\delta_{\vec{v}}^2})$ as follows:
\begin{align*}
&\|\displaystyle\sum_{i,l}2^{-(\vec{v}i+\vec{w}l)\cdot\vec{\mathbf{1}}}m_{i,l}(\xi)\ast f\|_{L^p_{\alpha}}\\
&\lesssim\displaystyle\displaystyle\sum_{k>0}\|\displaystyle\displaystyle\sum_{i}\sum_{l}2^{-\vec{v}\cdot\vec{\mathbf{1}}i-\vec{w}\cdot\vec{\mathbf{1}}l}m_{i,l}^{\vee}\ast P_{i,l,k}\ast f\|_{L^p_{\alpha}}+\displaystyle\sum_{k\leq0}\|\displaystyle\displaystyle\sum_{i}\sum_{l}2^{-\vec{v}\cdot\vec{\mathbf{1}}i-\vec{w}\cdot\vec{\mathbf{1}}l}m_{i,l}^{\vee}\ast P_{i,l,k}\ast f\|_{L^p_{\alpha}}\\
&\lesssim\displaystyle\displaystyle\sum_{k>0}\sum_{i}\displaystyle\displaystyle\sum_{l}2^{\alpha k}2^{-(1-\alpha\delta_{\vec{v}}^2)\vec{v}\cdot\vec{\mathbf{1}}i}2^{-(1-\alpha\delta_{\vec{w}}^2)\vec{w}\cdot\vec{1}l}\|m_{i,l}^{\vee}\ast \tilde{P}_{i,l,k}\ast f\|_{p}\\
&\qquad\qquad\qquad\qquad\qquad\qquad\qquad\qquad\quad\ \ +\displaystyle\sum_{k\leq0}\sum_{i}\displaystyle\displaystyle\sum_{l}2^{\alpha k}2^{-(1-\alpha\delta_{\vec{v}}^2)\vec{v}\cdot\vec{\mathbf{1}}i}2^{-(1-\alpha\delta_{\vec{w}}^2)\vec{w}\cdot\vec{1}l}\|m_{i,l}^{\vee}\ast \tilde{P}_{i,l,k}\ast f\|_{p}.
\end{align*}

Because $1-\alpha\delta_{\vec{v}}^2>0$ and $1-\alpha\delta_{\vec{w}}^2>0$, it follows that $\|m_{J}^{\vee}\ast f\|_{L^p_{\alpha}}\lesssim\|f\|_{L^p}$ for $min (\frac{1}{2\delta_{\vec{w}}^2},\frac{1}{2\delta_{\vec{v}}^2})<\frac{1}{p}\leq\frac{1}{2}$,
from the fact that $m_{i,l}(\xi)=O(|\frac{\xi_3}{2^{(\delta_{\vec{v}}^2\vec{\mathbf{1}}\cdot\vec{v}i+\delta_{\vec{w}}^2\vec{\mathbf{1}}\cdot\vec{w}l)}}|^{-1})$ and the same argument as in Proposition 1 for the case with $\delta_1 >\delta_2$. Furthermore, the interpolation between the above results and the fact that $\|\mathcal{A}_I\|_{BMO}\lesssim \|f\|_{\infty}$ yields the result for the range $\frac{1}{p}\leq min (\frac{1}{2\delta_{\vec{w}}^2},\frac{1}{2\delta_{\vec{v}}^2})$.
$\ $

\begin{rmk}

Let $ P(t_1,t_2)$ correspond to 
\begin{equation*}
\Lambda(P)=\{(m,0),(M,N)\}(m,M,N\geq 2)
\end{equation*}

$(1)$ $\delta_{\vec{v}}^1= \delta_{\vec{v}}^2, \delta_{\vec{w}}^1<\delta_{\vec{w}}^2$
 
 Assume that $ |\textrm{det}(\nabla^2_{t_1t_2} P(t_1,t_2))|\neq 0$. Then, $\mathcal{A}_J$ for $P(t_1,t_2)$ yields the same result as in Proposition 2 for the case with $\delta_{\vec{v}}^1< \delta_{\vec{v}}^2$ and $\delta_{\vec{w}}^1<\delta_{\vec{w}}^2$. 

$\ $

 $(2)$ $\delta_{\vec{v}}^1>\delta_{\vec{v}}^2, \delta_{\vec{w}}^1=\delta_{\vec{w}}^2$
 
 Assume that $ |\textrm{det}(\nabla^2_{t_1t_2} P(t_1,t_2))|\neq 0$. Then, $\mathcal{A}_J$ for $P(t_1,t_2)$ yields the same result as in Proposition 2 for the case with $\delta_{\vec{v}}^1> \delta_{\vec{v}}^2$ and $\delta_{\vec{w}}^1>\delta_{\vec{w}}^2$. 
 
 $\ $
 
In fact,
 
\begin{equation*}
\|\displaystyle\sum_{i,l}m_{J}^{\vee}\ast f\|_{L^p_{\alpha}}\lesssim\displaystyle\sum_{j\leq M }\|\displaystyle\sum_{l\geq M}m_{J}^{\vee}\ast f\|_{L^p_{\alpha}}+\displaystyle\sum_{l\leq M}\|\displaystyle\sum_{j\geq M}m_{J}^{\vee}\ast f\|_{L^p_{\alpha}}+\|\displaystyle\sum_{i,l\geq M}m_{J}^{\vee}\ast f\|_{L^p_{\alpha}}.
\end{equation*}

The first two terms are dealt with as in Proposition 1 and Remark 4. Thus, by the remark3 we can observe that the intersection of the $L^p\rightarrow L^p_{\alpha}$ ranges of the first two terms yields the desired result. In addition, we can deal with the last term using the same argument as in Proposition 2.

$\ $

$(3)$ The important implication of Proposition 2 for the case that $\delta_{\vec{v}}^1<\delta_{\vec{v}}^2$ and $\delta_{\vec{w}}^1<\delta_{\vec{w}}^2$ is that when $(M,N)\not\in\mathcal{R}$, the $(p,\alpha(p))$ ranges for $\mathcal{A}_J$ do not change even if we take an arbitrary vector $\vec{w}$.
 \end{rmk}

\section{Proof of the theorem}

Let $(\vec{n}_k)_{k}$ be the collection of normal vectors of $\mathcal{E}(P)$ and $\mathcal{E}(P^1)$, in ascending order of their slopes.   $$\mathbb{Z}^2_+=(\displaystyle\bigcup_{k}\{\vec{n}_ki|i\in \mathbb{Z}_+\})\cup(\displaystyle\bigcup_{k}\displaystyle\bigcup_{\vec{s}\in \mathcal{S}}\{\vec{s}+\vec{n}_ki+\vec{n}_{k+1}l| i,l\in \mathbb{Z}_+\})$$

For each decomposed index set, define $(\mathcal{A}_{I_k})_k$ and $(\mathcal{A}_{J_k^s})_{(s,k)}$ by
$$A_{I_k}=\sum_{(j_1,j_2)\in J_k}\int f(x_1-t_1, x_2-t_2, x_3-P(t_1,t_2))\eta(\frac{t_1}{2^{-j_1}}, \frac{t_2}{2^{-j_2}})dt, $$
$$A_{J_k^s}=\sum_{(j_1,j_2)\in J_k^s}\int f(x_1-t_1, x_2-t_2, x_3-P(t_1,t_2))\eta(\frac{t_1}{2^{-j_1}}, \frac{t_2}{2^{-j_2}})dt,$$
 where $J_k=\{\vec{n_k}i|i\in \mathbb{Z}_+\}$ and $J_k^s=\{s+\vec{n_k}i+\vec{n_{k+1}}l|i,l\in \mathbb{Z}_+\}$.

$\ $

Because $\mathcal{A}(f)=\displaystyle\sum_{k}\mathcal{A}_{I_k}(f)+\displaystyle\sum_{s}\displaystyle\sum_{k}\mathcal{A}_{J_k^s}(f)$, we can obtain the $L^p\rightarrow L^p_{\alpha}$ ranges for $\mathcal{A}$ by taking the intersection ranges of $L^p\rightarrow L^p_{\alpha}$ for all $\mathcal{A}_{I_k},\ \mathcal{A}_{J_k^s}$. By Lemma 1, we can find the first and second dominant monomials for each $\vec{n_k}$. Furthermore, we can observe that for fixed $k$, the $L^p\rightarrow L^p_{\alpha}$ ranges for $\mathcal{A}_{J_k^s}$ are the same for all $s$. This is because we can employ the cutoff function $\eta(\frac{\cdot}{2^{-s_1}},\frac{\cdot}{2^{-s_2}})$ in $m_J$ instead of $\eta(\cdot, \cdot).$ Thus, if we impose (2.1) and (2.2) on the polynomial $P(t_1,t_2)$, then we can apply Proposition 1 and Proposition 2 for each $\mathcal{A}_{I_k},\ \mathcal{A}_{J_k^s}$. Furthermore, by these propositions and the corresponding remarks, only taking the intersection of the $(\frac{1}{p},\alpha)$ ranges for $\mathcal{A}_{I_k}(f)$ will yield the results of Theorem 1. 

By the definitions of $\delta^{(m,n)}$ and $\mathcal{R}$, these can be defined for $n_s$. However, because we assume that min$(m_s,n_s)=m_s$, the set $\mathcal{R}$ is always the empty set for $n_s$, which easily follows from the definition of $\mathcal{R}$.

$\bullet$ $\mathcal{R}=\phi$

There exists $\vec{n}_l$ for which  $\delta_{\vec{n}_l}(t_1^{m_s})=\delta_{\vec{n}_l}(t_1^Mt_2^N)$, and by Proposition 1 it follows that
\begin{equation*}
\|\mathcal{A}_{I_l}(f)\|_{L_{\frac{2}{p}}^p}\leq\|f\|_{L^p}\ \textrm{if}\ \frac{1}{p}\leq\frac{1}{\delta_{\vec{n}_l}(t_1^{m_s})+\delta_{\vec{n}_l}(t_1^Mt_2^N)}
\end{equation*}
\begin{equation*}
\|\mathcal{A}_{I_l}(f)\|_{L_{\alpha}}\leq\|f\|_{L^p}\ \textrm{for}\ \alpha=\frac{1}{\delta_{\vec{n}_l}(t_1^{m_s})},\ \textrm{if}\ \frac{1}{\delta_{\vec{n}_l}(t_1^{m_s})+\delta_{\vec{n}_l}(t_1^Mt_2^N)}\leq\frac{1}{p}\leq\frac{1}{2}.
\end{equation*}

$\ $
 
 On the other hand, as in Lemma 4 (2), when $\delta_{\vec{v}_i}(t_1^{m_s})\leq\delta_{\vec{v}_i}(t_1^Mt_2^N)$, it holds that $(M,N)\not\in\mathcal{R}\ \textrm{if and only if}$
 \begin{align*}
  &\frac{1}{\delta_{\vec{v}_{i+1}}(t_1^{m_s})+\delta_{\vec{v}_{i+1}}(t_1^Mt_2^N)}>\frac{\delta_{\vec{v}_{i+1}}(t_1^Mt_2^N)-\delta_{\vec{v}_{i}}(t_1^Mt_2^N)}{\delta_{\vec{v}_i}(t_1^{m_s})\delta_{\vec{v}_{i+1}}(t_1^Mt_2^N)-\delta_{\vec{v}_i}(t_1^Mt_2^N)\delta_{\vec{v}_{i+1}}(t_1^{m_s})}\\
 (&\Leftrightarrow\delta_{\vec{v}_i}(t_1^{m_s})-\delta_{\vec{v}_{i+1}}(t_1^{m_s})> \delta_{\vec{v}_{i+1}}(t_1^Mt_2^N)-\delta_{\vec{v}_i}(t_1^Mt_2^N))\\
 (&\Leftrightarrow \frac{1}{\delta_{\vec{v}_{i+1}}(t_1^{m_s})+\delta_{\vec{v}_{i+1}}(t_1^Mt_2^N)}>\frac{1}{\delta_{\vec{v}_i}(t_1^{m_s})+\delta_{\vec{v}_i}(t_1^Mt_2^N)}).
 \end{align*}

Thus, comparing $p,\alpha$ ranges for large $p$, we can realize that for $i> l$, the $L^p\rightarrow L^p_{\alpha}$ ranges for $\mathcal{A}_{I_i}$ are wider than those for $\mathcal{A}_{I_l}$. As a result, we can obtain the desired the results by only comparing the $p,\alpha$ ranges for $\mathcal{A}_{I_k}$ associated with the normal vector $n_k$ of $\mathcal{E}(P)$.

$\bullet$ $\mathcal{R}\neq\phi$

As in Proposition 2 and Remark 5, the monomials $(M_1,N_1),..,(M_n,N_n)$ contained in $\mathcal{R}$ affect the $L^p\rightarrow L^p_{\alpha}$ ranges for $\mathcal{A}$. Thus, repeatedly applying Proposition 2 yields desired the results.

\begin{coro}
 Assume that $\frac{1}{\delta}<\frac{1}{2}$. $\mathcal{A}$ maps to $L^p\rightarrow L^q$ for the interior ranges of the $(\frac{1}{p},\frac{1}{q})$ diagrams, as in the following Figures.
\end{coro}

\begin{figure*}[ht]
  \includegraphics[width=0.5\textwidth]{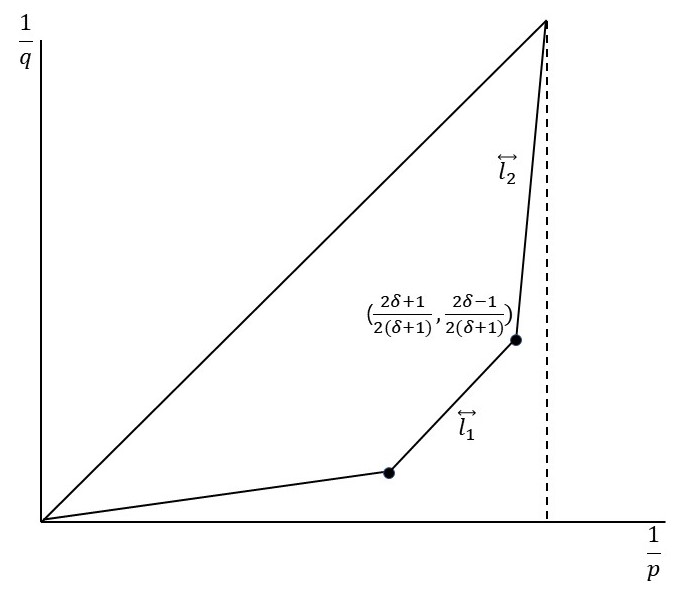}
\caption{$\mathcal{R}= \phi$}  
\label{5}
\end{figure*}



\begin{figure*}[ht]
  \includegraphics[width=0.7\textwidth]{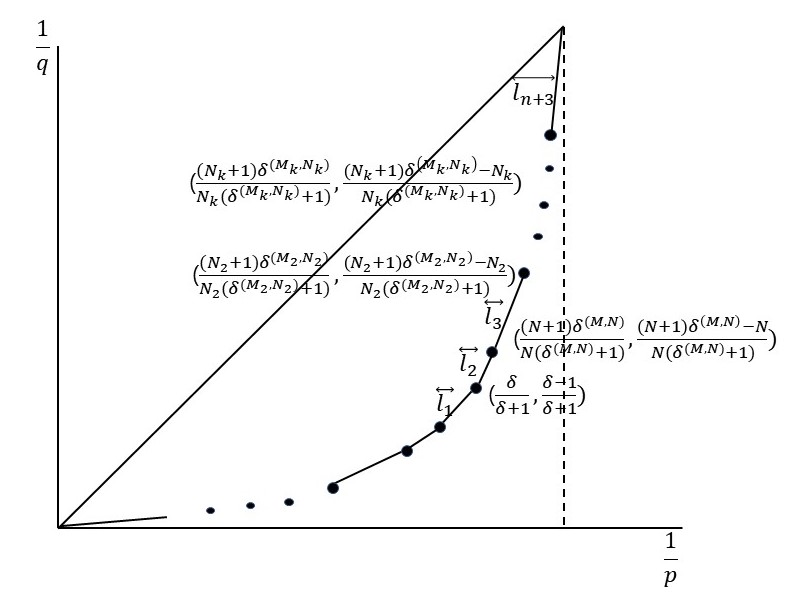}
\caption{$\mathcal{R}\neq \phi,$}       
\label{8}
\end{figure*}

\section{Necessary conditions}

We may assume that $\frac{1}{q}<1-\frac{1}{p}$, because otherwise we can utilize the duality of our operator.
We may assume that $\frac{1}{\delta}<\frac{1}{2}$, since the following necessary condition include the case $\frac{1}{\delta}\geq \frac{1}{2}$.

Necessary condition for $\overleftrightarrow{l}_1$:

We let
\begin{align*}
&Q_{\epsilon}=\{x : |x_1|<\epsilon^{\frac{N_1-N_2}{M_2N_1-M_1N_2}}, |x_2|< \epsilon^{\frac{M_2-M_1}{M_2N_1-M_1N_2}}, |x_3|<\epsilon\}\\& D_{\epsilon}=\{x : |x_1|<\epsilon^{\frac{N_1-N_2}{M_2N_1-M_1N_2}}, |x_2|< \epsilon^{\frac{M_2-M_1}{M_2N_1-M_1N_2}}, |x_3-P(x_1,x_2)|<\epsilon\}\\  &Q_{\epsilon}(x)=\{(t_1, t_2): |x_1-t_1|<\epsilon^{\frac{N_1-N_2}{M_2N_1-M_1N_2}}, |x_2-t_2|<\epsilon^{\frac{M_2-M_1}{M_2N_1-M_1N_2}},|x_3-P(t_1,t_2)|<\epsilon\}\\
&f=\chi_{Q_{\epsilon}},
\end{align*}

where $\overline{(M_1,N_1)(M_2,N_2)}$ is the edge intersecting $Y=X$.

Then, for each $x\in D_{\epsilon}$ and $t \in Q_{\epsilon}(x)$ it holds that
$$|P(x_1,x_2)-P(t_1,t_2)|\leq\displaystyle\sum_{(m,n)\in \Lambda(P)}|a_{mn}|(|(x_1^m||x_2^n-t_2^n)|+|x_1^m-t_1^m||t_2^n|)\lesssim \epsilon.$$

Indeed, 
\begin{equation*}
|(x_1^m||x_2^n-t_2^n)|+|x_1^m-t_1^m||t_2^n| \lesssim \epsilon^{\frac{N_1-N_2}{M_2N_1-M_1N_2}m+\frac{M_2-M_1}{M_2N_1-M_1N_2}n}=\epsilon^{\frac{(N_1-N_2)(m-M_1)+(M_2-M_1)(n-N_1)}{M_2N_1-M_1N_2}+1}.
\end{equation*}
Thus, if $(m,n)(\in \Lambda(P))$ lies on $\overline{(M_1,N_1)(M_2,N_2)}$, then \begin{equation*}
\frac{(N_1-N_2)(m-M_1)+(M_2-M_1)(n-N_1)}{M_2N_1-M_1N_2}+1=1.
\end{equation*}
Otherwise, it is greater than $1$. 

$\ $

Therefore, for each $x\in D_{\epsilon}$,
\begin{equation*}
\mathcal{A}(f)=\int \chi_{Q_{\epsilon}}(x_1-t_1,x_2-t_2,x_3-P(t_1,t_2))\psi(t)dt \geq\int_{Q_{\epsilon}(x)}dt
= |Q_{\epsilon}(x)|=\epsilon^{\frac{N_1-N_2+M_2-M_1}{M_2N_1-M_1N_2}},
\end{equation*}
and if $\mathcal{A}$ maps $L^p$ into $L^q$ then $(\displaystyle\int_{D_{\epsilon}}|\mathcal{A}(f)(x)|^qdx)^{\frac{1}{q}}\lesssim\|f\|_{p}$. Also, it holds that
$$|\epsilon^{\frac{N_1-N_2+M_2-M_1}{M_2N_1-M_1N_2}+1}|^{\frac{1}{p}}\gtrsim|Q_{\epsilon}|^{\frac{1}{p}}=\|f\|_p$$
\begin{equation*}
(\int_{D_{\epsilon}}|\mathcal{A}(f)(x)|^qdx)^{\frac{1}{q}}\gtrsim \epsilon^{\frac{N_1-N_2+M_2-M_1}{M_2N_1-M_1N_2}}|D_{\epsilon}|^{\frac{1}{q}}=\epsilon^{\frac{N_1-N_2+M_2-M_1}{M_2N_1-M_1N_2}}|\epsilon^{\frac{N_1-N_2+M_2-M_1}{M_2N_1-M_1N_2}+1}|^{\frac{1}{q}}.
\end{equation*}
Letting $\epsilon\rightarrow0$ and comparing the exponents, it is necessary that $$\frac{1}{q}\gtrsim \frac{1}{p}-\frac{N_1-N_2+M_2-M_1}{M_2N_1-M_1N_2+N_1-N_2+M_2-M_1}.$$
Necessary condition for $\overleftrightarrow{l}_2$:
We let
\begin{align*}
&Q_{\epsilon}=\{x : |x_1|<\epsilon^{1-\frac{m_s-1}{M}}, |x_2|< 1, |x_3|<\epsilon\}\\ &D_{\epsilon}=\{x : |x_1|<\epsilon^{\frac{1}{M}}, |x_2|< 1, |x_3-P(x_1,x_2)|<\epsilon\}\\  &Q_{\epsilon}(x)=\{(t_1, t_2): |x_1-t_1|<\epsilon^{1-\frac{m_s-1}{M}}, |x_2-t_2|<1, |x_3-P(t_1,t_2)|<\epsilon\}\\
&f=\chi_{Q_{\epsilon}},
\end{align*}

where $m_s$ and $M$ are as in Definitions 2 and 3.
Then, for each $x\in D_{\epsilon}$ and $t \in Q_{\epsilon}(x)$,
$$|P(x_1,x_2)-P(t_1,t_2)|\leq\displaystyle\sum_{(m,n)\in \Lambda(P)}|a_{mn}|(|(x_1^m||x_2^n-t_2^n)|+|x_1^m-t_1^m||t_2^n|)\lesssim \epsilon$$

Therefore, for each $x\in D_{\epsilon}$, it holds that

$\mathcal{A}(f)=\displaystyle\int \chi_{Q_{\epsilon}}(x_1-t_1,x_2-t_2,x_3-P(t_1,t_2))\psi(t)dt \geq\displaystyle\int_{Q_{\epsilon}(x)}1dt =\epsilon^{1-\frac{m_s-1}{M}},$
and if $\mathcal{A}$ maps $L^p$ into $L^q$, then $(\displaystyle\int_{D_{\epsilon}}|\mathcal{A}(f)(x)|^qdx)^{\frac{1}{q}}\lesssim\|f\|_{p}$ and

$$\|f\|_p=|Q_{\epsilon}|^{\frac{1}{p}}\lesssim|\epsilon^{2-\frac{m_s-1}{M}}|^{\frac{1}{p}}$$
$$(\int_{D_{\epsilon}}|\mathcal{A}(f)(x)|^qdx)^{\frac{1}{q}}\gtrsim \epsilon^{1-\frac{m_s-1}{M}}|D_{\epsilon}|^{\frac{1}{q}}\gtrsim \epsilon^{1-\frac{m_s-1}{M}}\epsilon^{(1+\frac{1}{M})\frac{1}{q}}.$$
Letting $\epsilon\rightarrow0$ and comparing the exponents, it is necessary that

$$(1+\frac{1}{M})\frac{1}{q}\geq(2-\frac{m_s-1}{M})\frac{1}{p}-(1-\frac{m_s-1}{M}).$$
Necessary condition for $\overleftrightarrow{l}_{k}$ ($k\geq3$):

We let
\begin{align*}
&Q_{\epsilon}=\{x : |x_1|<\epsilon^{1-(m_s-1)\frac{N_k-N_{k+1}}{M_{k+1}N_k-M_kN_{k+1}}}, |x_2|< \epsilon^{\frac{M_{k+1}-M_k}{M_{k+1}N_k-M_kN_{k+1}}}, |x_3|<\epsilon\}\\ &D_{\epsilon}=\{x : |x_1|<\epsilon^{\frac{N_k-N_{k+1}}{M_{k+1}N_k-M_kN_{k+1}}}, |x_2|< \epsilon^{\frac{M_{k+1}-M_k}{M_{k+1}N_k-M_kN_{k+1}}},|x_3-P(x_1,x_2)|<\epsilon\}\\  
&Q_{\epsilon}(x)=\{(t_1, t_2): |x_1-t_1|<\epsilon^{1-(m_s-1)\frac{N_k-N_{k+1}}{M_{k+1}N_k-M_kN_{k+1}}},\\
&\qquad\qquad\qquad\qquad\qquad\qquad|x_2-t_2|<\epsilon^{\frac{M_{k+1}-M_k}{M_{k+1}N_k-M_kN_{k+1}}},
|x_3-P(t_1,t_2)|<\epsilon\}\\
&f=\chi_{Q_{\epsilon}}, 
\end{align*}
where the $(M_k,N_k)$ terms are enumerated as in Definition 4.

Then, for each $x\in D_{\epsilon}$ and $t \in Q_{\epsilon}(x)$ it holds that
$$|P(x_1,x_2)-P(t_1,t_2)|\leq\displaystyle\sum_{(m,n)\in \Lambda(P)}|a_{mn}|(|(x_1^m||x_2^n-t_2^n)|+|x_1^m-t_1^m||t_2^n|)\lesssim \epsilon.$$
Indeed, by comparing the slopes of $\overline{(m_s,0)(M_k,N_k)}$ and $\overline{(m_s,0)(M_{k+1},N_{k+1})}$, we can observe that 
\begin{equation*}
    (m_s-M_{k+1})N_k<(m_s-M_k)N_{k+1},
\end{equation*}
which implies that $\frac{N_k-N_{k+1}}{M_{k+1}N_k-M_kN_{k+1}}<1-(m_s-1)\frac{N_k-N_{k+1}}{M_{k+1}N_k-M_kN_{k+1}}$. Thus, if $x\in D_{\epsilon}$ and $t\in Q_{\epsilon}(x)$, then $|t_1|\lesssim \epsilon^{\frac{N_k-N_{k+1}}{M_{k+1}N_k-M_kN_{k+1}}}$. Thus, $|x_1^{m_s}-t_1^{m_s}|\lesssim \epsilon$, and for $m,n \geq 2$
\begin{multline*}
|x_1^m||x_2^n-t_2^n|+|x_1^m-t_1^m||t_2^n|\lesssim \epsilon^{m\frac{N_k-N_{k+1}}{M_{k+1}N_k-M_kN_{k+1}}+n\frac{M_{k+1}-M_k}{M_{k+1}N_k-M_kN_{k+1}}}\\+\epsilon^{1-(m_s-1)\frac{N_k-N_{k+1}}{M_{k+1}N_k-M_kN_{k+1}}+(m-1)\frac{N_k-N_{k+1}}{M_{k+1}N_k-M_kN_{k+1}}+n\frac{M_{k+1}-M_k}{M_{k+1}N_k-M_kN_{k+1}}}
\end{multline*}

Comparing the slopes of $\overline{(M_{k+1}N_{k+1})(M_k,N_k)}$ and $\overline{(m,n)(M_k,N_k)}$, the first term is observed to be less than $\epsilon$. Comparing the slopes of 
$\overline{(m_s,0)(m,n)}$ and $\overline{(M_k,N_k)(M_{k+1}N_{k+1})}$, we can observe that
\begin{equation*}
    (m-m_s)(N_k-N_{k+1})\geq n(M_k-M_{k+1}),
\end{equation*}
which implies that
\begin{equation*}
    1+(m-m_s)\frac{N_k-N_{k+1}}{M_{k+1}N_k-M_kN_{k+1}}+n\frac{M_{k+1}-M_k}{M_{k+1}N_k-M_kN_{k+1}}\geq 1.
\end{equation*}
Thus, 
\begin{equation*}
    \epsilon^{1-(m_s-1)\frac{N_k-N_{k+1}}{M_{k+1}N_k-M_kN_{k+1}}+(m-1)\frac{N_k-N_{k+1}}{M_{k+1}N_k-M_kN_{k+1}}+n\frac{M_{k+1}-M_k}{M_{k+1}N_k-M_kN_{k+1}}}\leq\epsilon.
\end{equation*}
Therefore, for each $x\in D_{\epsilon},$

\begin{multline*}
\mathcal{A}(f)=\displaystyle\int \chi_{Q_{\epsilon}}(x_1-t_1,x_2-t_2,x_3-P(t_1,t_2))\psi(t)dt \geq\displaystyle\int_{Q_{\epsilon}(x)}dt = |Q_{\epsilon}(x)|\\=\epsilon^{1-(m_s-1)\frac{N_k-N_{k+1}}{M_{k+1}N_k-M_kN_{k+1}}+\frac{M_{k+1}-M_k}{M_{k+1}N_k-M_kN_{k+1}}},
\end{multline*}
and if $\mathcal{A}$ maps $L^p$ into $L^q$, then $(\displaystyle\int_{D_{\epsilon}}|\mathcal{A}(f)(x)|^qdx)^{\frac{1}{q}}\lesssim\|f\|_{p}$ and

$$\|f\|_p=|Q_{\epsilon}|^{\frac{1}{p}}\lesssim|\epsilon^{2-(m_s-1)\frac{N_k-N_{k+1}}{M_{k+1}N_k-M_kN_{k+1}}+\frac{M_{k+1}-M_k}{M_{k+1}N_k-M_kN_{k+1}}}$$
$$\epsilon^{1-(m_s-1)\frac{N_k-N_{k+1}}{M_{k+1}N_k-M_kN_{k+1}}+\frac{M_{k+1}-M_k}{M_{k+1}N_k-M_kN_{k+1}}}|D_{\epsilon}|^{\frac{1}{q}}\lesssim(\int_{D_{\epsilon}}|\mathcal{A}(f)(x)|^qdx)^{\frac{1}{q}}.$$
Letting $\epsilon\rightarrow0$ and comparing the exponents, it is necessary that
\begin{align*}
&(2-(m_s-1)\frac{N_k-N_{k+1}}{M_{k+1}N_k-M_kN_{k+1}}+\frac{M_{k+1}-M_k}{M_{k+1}N_k-M_kN_{k+1}})\frac{1}{p}
\\&\leq(1-(m_s-1)\frac{N_k-N_{k+1}}{M_{k+1}N_k-M_kN_{k+1}}+\frac{M_{k+1}-M_k}{M_{k+1}N_k-M_kN_{k+1}})\\&\qquad\qquad\quad\qquad\qquad\qquad\qquad\qquad\qquad\qquad+\frac{1}{q}(1+\frac{N_k-N_{k+1}+M_{k+1}-M_k}{M_{k+1}N_k-M_kN_{k+1}}).\ \ (5.1)
\end{align*}

Now, we check that this is a necessary condition for $\overleftrightarrow{l}_{k+1}$. Because $\overleftrightarrow{l}_{k+1}$ passes through 
\begin{equation*}
    (\frac{(N_k+1)\delta^{(M_k,N_k)}}{N_k(\delta^{(M_k,N_k)}+1)}, \frac{(N_k+1)\delta^{(M_k,N_k)}-N_k}{N_k(\delta^{(M_k,N_k)}+1)})
\end{equation*}
\begin{equation*}
(\frac{(N_{k+1}+1)\delta^{(M_{k+1},N_{k+1})}}{N_{k+1}(\delta^{(M_{k+1},N_{k+1})}+1)}, \frac{(N_{k+1}+1)\delta^{(M_{k+1},N_{k+1})}-N_{k+1}}{N_{k+1}(\delta^{(M_{k+1},N_{k+1})}+1)}),
\end{equation*}
we insert these points into (5.1).

Let $\frac{N_k-N_{k+1}}{M_{k+1}N_k-M_kN_{k+1}}=A$ and $\frac{M_{k+1}-M_k}{M_{k+1}N_k-M_kN_{k+1}}=B$, for notational simplicity. After some calculations, we obtain the following:
\begin{align*}
    &(2-(m_s-1)A+B)\frac{\delta^{(M_k,N_k)}-N_k}{N_k(\delta^{(M_k,N_k)}+1)}=A+B+\frac{\delta^{(M_k,N_k)}-2N_k}{N_k(\delta^{(M_k,N_k)}+1)}(1+A+B)\\
    &\Leftrightarrow(1-m_sA)\frac{\delta^{(M_k,N_k)}-N_k}{N_k(\delta^{(M_k,N_k)}+1)}=(A+B)-\frac{1}{(\delta^{(M_k,N_k)}+1)}(1+A+B)\\
    &\Leftrightarrow(1-m_sA)(\delta^{(M_k,N_k)}-N_k)=N_k\delta^{(M_k,N_k)}(A+B)-N_k\\
    &\Leftrightarrow(1-m_sA)\delta^{(M_k,N_k)}=N_kA(\delta^{(M_k,N_k)}-m_s)+N_k\delta^{(M_k,N_k)}B\\
    &\Leftrightarrow(AM_k-m_sA)\delta^{(M_k,N_k)}=N_kA(\delta^{(M_k,N_k)}-m_s)
\end{align*}
We used the fact $AM_k+BN_k=1$ for the final line. Finally, let us insert 
\begin{equation*}
    \delta^{(M_k,N_k)}=\frac{N_km_s}{N_k-M_k+m_s}
\end{equation*}
into the last equality. 

Finally, for small $\epsilon$ we let

$\left\{\begin{array}{l} Q_{\epsilon}=\{x : |x_1|<\epsilon, |x_2|< \epsilon, |x_3|\lesssim\epsilon\}\\ D_{\epsilon}=\{x: |x_1|<1, |x_2|<1, |x_3-P(x_1,x_2)|<\epsilon\}\\  Q_{\epsilon}(x)=\{(t_1, t_2): |x_1-t_1|<\epsilon, |x_2-t_2|<\epsilon, |x_3-P(t_1,t_2)|\lesssim\epsilon\}\\
f=\chi_{Q_{\epsilon}} \end{array}\right.$ 

Then, for each $x\in D_{\epsilon}$ and $t \in Q_{\epsilon}(x)$,
$$|P(x_1,x_2)-P(t_1,t_2)|\leq\displaystyle\sum_{(m,n)\in \Lambda(P)}|a_{mn}|(|(x_1^m||x_2^n-t_2^n)|+|x_1^m-t_1^m||t_2^n|)\lesssim \epsilon,$$

Furthermore, for each $x\in D_{\epsilon}$, 

$$A(f)=\int \chi_{Q_{\epsilon}}(x_1-t_1,x_2-t_2,x_3-P(t_1,t_2))\psi(t)dt \gtrsim \displaystyle\int_{Q_{\epsilon}(x)}dt=|Q_{\epsilon}(x)|=\epsilon^2,$$

and if $\mathcal{A}$ maps $L^p$ into $L^q$, then 

$$\epsilon^{\frac{3}{p}}\gtrsim |Q_{\epsilon}|^{\frac{1}{p}}=\|f\|_p\gtrsim(\int_{D_{\epsilon}}|A(f)(x)|^qdx)^{\frac{1}{q}}\gtrsim \epsilon^2|D_{\epsilon}|^{\frac{1}{q}}=\epsilon^{2+\frac{1}{q}}.$$
Letting $\epsilon\rightarrow 0$ and comparing the exponents, it is necessary that $\frac{1}{q}\geq\frac{3}{p}-2.$


%
%



\end{document}